\documentclass[12pt,a4paper]{amsart}
\usepackage{amssymb}
\usepackage{amsmath}
\usepackage{latexsym}

\usepackage[T1]{fontenc}
\usepackage[utf8]{inputenc}
\usepackage[german,english]{babel}
\usepackage{soul}

\usepackage[babel,german=guillemets]{csquotes}
\usepackage[backend=bibtex8,style=authoryear]{biblatex}
\bibliography{auffassungen.bib}
\usepackage{lmodern}

\usepackage[bottom]{footmisc}

\begin{document}

\pagestyle{plain}

\author{Hans Joachim Burscheid $\cdot$ Horst Struve}
\address{Universität zu Köln, Institut für Mathematikdidaktik \\ Gronewaldstr.2, 50931 Köln \\ Tel.: 0049--221--470--4750, Fax.: 0049-221-470-4985 \\ e-mail: h.burscheid@uni-koeln.de, h.struve@uni-koeln.de}

\title{
\large{Mathematisches Wissen,} \\ 
\vspace{3mm}
\small{erworben in empirischen Theorien}}

\selectlanguage{german}
\begin{abstract}
Im Anschluß an die Ausarbeitung einzelner Themen der schulischen Elementarmathematik als Inhalte empirischer Theorien wird in diesem Beitrag die Frage aufgegriffen, ob die damit verbundene Auffassung von Mathematik sich unter etablierten Auffassungen wiederfindet, und wie sie sich stützen läßt. Eine lernpsychologische Position und eine Unterrichtskonzeption werden herangezogen, die zu der hier entwickelten Auffassung von Mathematik und der über ihr Erlernen passen. 
\end{abstract}

\maketitle
\vspace{-1cm}

\selectlanguage{english}
\begin{abstract}
Following the processing of individual topics of elementary school mathematics as content of empirical theories the question is adressed wether the associated conception of mathematics finds itself under established concepts, and how it can be supported. A position of learning psychology and a teaching concept will be used that fit in the position of mathematics and of the learning of mathematics developed here. 
\end{abstract}

\vspace{5mm}
\thispagestyle{empty}

\noindent \sloppy In [2009] haben wir den Standpunkt ausgearbeitet, daß die schulübliche Elementarmathematik im Rahmen empirischer Theorien erworben werde und nicht in der Form mathematischer Theorieteile. \emph{Empirische Theorie} wird dabei als Oberbegriff aller solcher Theorien verstanden, die ein Stück WELT beschreiben. Paradigmatische Beispiele derartiger Theorien sind die Theorien der Naturwissenschaften. Zur Begründung dieser Auffassung wurde darauf verwiesen, daß die Mathematik, die Schüler\footnote{und natürlich auch Schülerinnen} erlernen, ihnen dazu dient, Probleme zu bewältigen, die ihre Umwelt ihnen stellt. Diese Probleme haben durchweg empirischen Charakter, sodaß die Begriffe, die die Schüler erwerben, eine starke ontologische Bindung haben, wie Hans Freudenthal es ausdrückte. Diese Auffassung über das Erlernen der Elementarmathematik wird z.B. auch von Heinz Griesel geteilt [2013]. Daß bestimmtes kindliches Wissen beschrieben werden kann, als sei es in Theorien organisiert, ist eine Einsicht, die Psychologen bei der Beobachtung von Kleinkindern erhielten (Gopnik und Meltzoff [1997]), die aber auch schon früher bei Schülern gewonnen wurde (Resnick [1983]). \smallskip

Wozu wir uns bislang nicht geäußert haben, ist die Frage, ob das im Rahmen empirischer Theorien erworbene mathematische Wissen sich einer Auffassung von Mathematik zuordnen läßt, die als \enquote{etabliert} bezeichnet werden kann. Auch lernpsychologische Fragen oder die Frage nach einer geeigneten Unterrichtskonzeption haben wir bislang nicht angesprochen. Diesen Fragen werden wir uns in diesem Beitragc widmen. \medskip

\noindent \emph{Bemerkung}: Im folgenden ist es unvermeidlich, in ausgiebiger Weise Zitate heranzuziehen, da sich nur so belegen läßt, daß die genannten Autorinnen und Autoren die ihnen zugewiesenen Auffassungen auch vertreten haben.\medskip

Auf dem \emph{12$^{th}$ International Congress on the History of Science} im Jahre 1968 in Paris wurde eine Initiative angestoßen, die die Gründung der \emph{International Commission on the History of Mathematics} zur Folge hatte. In Band 2 der von der Kommission herausgegebenen Zeitschrift \emph{Historia Mathematica} weist der erste Vorsitzende der Kommission, Kenneth O. May, auf Fehler und Mängel bisheriger mathematischer Geschichtsschreibung hin. Dort heißt es:

\begin{small}
\enquote{These vices arise, as do most of those I shall discuss, from looking at mathematics as a timeless, static, logical structure gradually revealed to man, instead of a historically evolving human phenomenon.} [1975, p. 187]
\end{small} 

Die Kommission führte 1974 einen \emph{Workshop on the Evolution of Modern Mathematics} durch, dessen Bericht ebenfalls in Band 2 der \emph{Historia Mathematica} festgehalten ist. In einem schriftlichen Beitrag zu diesem Workshop argumentierten Nancy Kopell und Gabriel Stolzenberg:

\begin{small}
\enquote{The arbitrary limits to which we refer have their origins in the revolutionary change in mathematical world--view which took place during the
period 1870 -- 1930 (i.e., the period that begins with Weierstrass’ \enquote{rigorous} codification of analysis and ends with Hilbert’s formalist program for \enquote{saving} it). Before 1870, pure mathematics had a flourishing \emph{empirico--inductive} tradition which included, though not as its \enquote{foundations}, a growing \emph{logico--deductive} component. The hallmark of the empirico--inductive tradition was its primary concern with the phenomena of mathematics. Its theories where theories about the phenomena, just as in physical theory. (Gauss’ \emph{Disquisitiones Arithmeticae} is a classic work in this tradition) By 1930, the situation has changed completely. By then, there was a nearly universal acceptance of the modern strictly logico–deductive conception of pure mathematics and the empirico–inductive tradition had been virtually suppressed. In empirico–inductive mathematics, the meaning of any theorem is in what it tells us about the phenomana it describes. By contrast, modern logico–deductive mathematics is marked by a detachment of the theorems from the phenomena (which are now demoted to the status of \enquote{examples} or \enquote{illustrations} of the theorems).} [1975, p. 519]
\end{small}

 Auf dem Internationalen Mathematikerkongreß im Jahre 1950 in Cambridge MA formulierte der Mathematikphilosoph Raymond Wilder in seinem Vortrag:

\begin{small}
\enquote{As a body of knowledge mathematics is ..... a part of our culture, our \emph{collective} possession. [1986, p. 188]

A culture is the collection of customs, rituals, beliefs, tools, mores, etc., which we may call \emph{cultural elements}, possessed by a group of people, such as a primitive tribe or the people of North America. Generally it is not a fixed thing but changing with the course of time, forming what can be called a \enquote{culture stream}. It is handed down from one generation to another, constituting a seemingly living body of tradition .....} [p. 187]
\end{small}

Hinsichtlich der Bemühungen um eine \enquote{Grundlagentheorie} der Mathematik, wie sie in den ersten Jahrzehnten des 20sten Jahrhunderts angestrengt wurden, führte er aus:

\begin{small}           
\enquote{If the culture concept tells us anything, it should teach us that the first rule for setting up any Foundation theory is that it should only attempt to encompass specific portions of the field as it is known in our culture. [p. 195]

Problems of mathematical existence, for example, can never be settled by appeal to any mathematical dogma. Indeed, they have no validity except as related to special foundations theories. (...) \emph{Because of its cultural basis, there is no such thing as the absolute in mathematics; there is only the relative.}} [p. 196]
\end{small}

Mathematik als kulturelle Errungenschaft zu verstehen betont auch Morris Kline. Er schreibt:

\begin{small}
\enquote{It would appear as though mathematics is the creation of human, fallible minds rather than a fixed, eternally existing body of knowledge. The subject seems very much dependent upon the creator. As Alfred North Whitehead put it, \enquote{The science of pure mathematics may claim to be the most original creation of the human spirit.} Only the \emph{relatively} universal acceptance of mathematics (as opposed to the acceptance of religions, political, and ethical doctrines) may lure us into granting that subject an objective existence.} [1962, p. 665] 
\end{small}

Eine Sichtweise der Mathematik, die gestützt wurde durch ihre herausragende Rolle  in Produktion, Technik und dem gesamten sozio--ökonomischen Prozeß, wie sie der Mathematikhistoriker Dirk Jan Struik (1894 -- 2000) --- offenbar ein überzeugter Vertreter des dialektischen Materialismus --- immer wieder betonte, der seit 1926 am MIT lehrte [1936, 1942]. \vspace{3mm}

In den 1960er und vor allem in den 1970er Jahren ergriff die Diskussion über die philosophische Sicht der Mathematik weite Kreise der Mathematiker. Christopher Ormell gibt dieser Entwicklung eine Bedeutung vergleichbar der Grundlagenkrise zu Beginn des 20. Jahrhunderts [1992, pp. 5/6]. Aber nicht nur die eher pragmatisch orientierte Sicht der angelsächsischen Mathematiker trug dazu bei. Ein wesentlicher
Anstoß ging von dem ungarischen Mathematiker und Wissenschaftsphilosophen Imre Lakatos (1922 -- 1974) und dem amerikanischen Philosophen Hilary Putnam aus. 

Lakatos unterscheidet \emph{euklidische} Theorien --- solche, in denen endlich viele Wahrheitssetzungen (die Axiome) an die Spitze gestellt werden, aus denen man die Sätze der Theorie nach logischen Regeln ableitet --- von \emph{nicht--euklidischen} Theorien. Diese nennt er \emph{quasi--empirisch}. Während er einräumt, daß \enquote{lokale} mathematische Theorien wie z.B. die Gruppentheorie dem Anspruch genügen, euklidisch zu sein, bestreitet er dies für \enquote{umfassende} Theorien, was er mit dem Scheitern des Logizismus und des Formalismus begründet, der beiden großen Versuche, die Probleme zu lösen, die sich durch das Auftauchen der Paradoxien der Mengenlehre aufgetan hatten. Zu letzterem Punkt schreibt Nicholas D. Goodman:

\begin{small}
\enquote{The independence results, the proliferation of large cardinal axioms, and
the construction of increasingly bizarre models for set theory have made
mathematicians realize how weak their set--theoretic intuition actually is.}
[1979, p. 549]
\end{small}

Als quasi--empirische Wissenschaft sieht Lakatos die Mathematik wie eine Naturwissenschaft [1982, 2. Kap.]. Auch Philip J. Davis argumentiert:

\begin{small}
\enquote{The peculiar arguments made here ..... lead to the conclusion that mathematics, in some of his aspects, takes on the nature of an experimental
science.} [1972, p. 253] 
\end{small}

Diese Sichtweise wird durch Argumente gestützt, die sich auf die Entstehung der Mathematik beziehen. So heißt es bei Warwick Sawyer:
\begin{small}
\enquote{For mathematics began with the study of real things, arithmetic for
counting them, geometry and trigonometry for measuring them, ..... [1992,
p. 68]

\noindent Spengler has pointed out, I believe correctly, that the Greeks had no
concept of empty space. Their ’geometry’ was mathematical physics, the study of the shapes and sizes of material objects ..... The reasoning was in terms of familiar objects, and many properties were taken for granted without being analysed. Socrates’ demonstration of a particular case of Pythagoras’ Theorem, which is highly convincing on the level of common sense, assumes all kinds of things about the nature of area and the way in which squares fit together. We find that reasoning at this level is very effective, and Euclid’s geometry has for centuries provided a basis for surveying and manufacturing design ..... completely reliable to the degree of accuracy required. [pp. 69/70]

\noindent Our belief is that a result such as 2 + 3 = 5 obtained by experiment with
pebbles or blocks will apply equally well to any kind of objects that do not fuse, split or evaporate. ..... arithmetic definitely tells us experimentally verifiable facts about the world.} [p. 70] 
\end{small}

\noindent Nahezu gleichlautend argumentiert Kline:

\begin{small}
\enquote{The basic concepts of the elementary branches of mathematics are abstractions from experience: Whole numbers and fractions were certainly suggested by obvious physical counterparts. [1962, pp. 660/661]

Fortunately, every abstraction is ultimately derived from, and therefore understandable in terms of, intuitively meaningful objects or phenomena. The mind does play its part in the creation of mathematical concepts, but the mind does not function independently of the outside world. Indeed the mathematician who treats concepts that have no physically real or intuitive origins is almost surely talking nonsense.} [p. 31] 
\end{small}

\noindent Saunders Mac Lane formuliert wie folgt:

\begin{small}
\enquote{I assert that subjects of Mathematics are extracted from the environment;
that is, from activities, phenomena, or science --- and that they are then later applied to that -- or other --- environments. Thus number theory is
\enquote{extracted} from the activity of counting, and geometry is extracted from
motion and shaping. (...) I have deliberately chosen this work (word; die Verf.) \enquote{extraction} to be close to the more familiar word \enquote{abstraction} --- and with the intent that the Mathematical subject resulting from an extraction
is indeed abstract. Mathematics is not \enquote{about} human activity, phenomena, or science. It is about the extractions and formalisations of ideas --- and their manifold consequences.} [1986, p. 418] \vspace{3mm}
\end{small}

Welche Ausprägung der Idee, mathematische Objekte oder Zusammenhänge aus menschlichen Aktivitäten herauszufiltern man bevorzugt --- bei Davis heißt es -- ähnlich wie bei Kline 

\begin{small}
\enquote{..... geometry can be regarded as an abstraction, a distillation, a formalization, an intellectualization of the visual and kinesthetic experience of
space}, [1994, p. 166]
\end{small}

und welche Tätigkeiten auch zugrunde liegen, die implizite Regelmäßigkeiten enthalten, es ist sorgfältig darauf zu achten, daß man nicht der Gefahr erliegt,

\begin{small}
\enquote{..... that one will see only the context and miss the mathematics.}
(Thomas [1996, p. 13])
\end{small}

Denn Mathematik handelt nicht nur von den Mustern, die aus dem realen Handeln oder den realen Objekten herausgefiltert werden, sondern impliziert deren Rückbindung an die Realität, denn nur diese

\begin{small}
\enquote{..... justifies Frege’s dictum that it is applicability alone that raises arithmetic from the rank of a game to that of a science.} (Dummett [1994, p. 14])
\end{small}

Oder mit den Worten von Robert Thomas:

\begin{small}
\enquote{In my opinion mathematics is not just about these extractions and their
formulization, but about the reattachment of these ideas to \enquote{human activity, phenomena, or science}. The whole enterprise involves the detachments of
insights, the study of the insights, and the re--attachment of the insights, or as it normally called the application of the mathematics.} [1996, p. 16]
\end{small}

Freudenthal formuliert dies sehr einprägsam:

\begin{small}
\enquote{Das --- ich meine die Wirklichkeit --- ist das Skelett, an dem die
Mathematik sich festsetzt, .....} [1973, S. 77] \medskip
\end{small}

Lakatos’ einflußreichstes Werk war \enquote{Beweise und Widerlegungen} (1976 posthum im englischen Original erschienen). Obwohl er ausschließlich aus wissenschaftsphilosophischer Sicht argumentierte, entfachte er -- möglicherweise durch Äußerungen wie der folgenden

\begin{small}
--- \enquote{Es braucht mehr als die Antinomien (der Mengenlehre, die Verf.) und die Gödelschen Ergebnisse, um die Philosophen dazu zu bringen, die empirischen Seiten der Mathematik ernst zu nehmen und eine Philosophie des kritischen Fallibilismus zu entwickeln, die sich nicht von den sogenannten Grundlagen inspirieren läßt, sondern von der Entwicklung der mathematischen Erkenntnis.} [1982, S. 41] ---
\end{small}

unter den Mathematikern eine Diskussion über eine philosophische Sicht der Mathematik, die sich enger an die Tätigkeit des Mathematikers anlehnt. Erstmalig wurde dabei auch die mathematische Lehre ins Spiel gebracht, sogar darauf verwiesen, daß während der 1960er Jahre der Versuch, Axiomatik und mengentheoretische Schreibweise in der Schule einzuführen, Ausdruck der Auffassung war, daß die gesamte Mathematik aus axiomatisch formulierten Lehrgebäuden in mengentheoretischer Sprache bestehe (Hersh [1979, p. 33]). Reuben Hersh plädierte nachhaltig dafür, daß die Philosophie der Mathematik unter Zurückweisung jeder formalistischen oder platonischen Sichtweise das mathematische Wissen wie jedes andere menschliche Wissen als fehlbar (fallible), verbesserbar (corrigible), vorläufig (tentative) und entwickelbar(evolving) darstellen soll. [p. 43] \vspace{4mm} 

Neben diese überwiegend von in der Forschung tätigen Mathematikern vertretene Auffassung von Mathematik stellen wir eine zweite, die Mathematik aus der Sicht eines Wissenschaftsphilosophen beleuchtet. Sie wurde von dem Engländer Philip Kitcher vorgelegt, dessen  mathematik--philosophische Position die sich wandelnde Auffassung von Mathematik noch einmal zuspitzte [1984]. Kitcher war durch die Beschäftigung mit den Gedanken des englischen Philosophen John Stuart Mill (1806 -- 1873) angeregt worden, mathematisches Wissen enger an  die Realität zu binden, als er untersuchte, in welchem Sinne die Auffassung Mills zu verstehen sei, arithmetisches Wissen sei empirisch.

\begin{small}
\enquote{He (Mill, die Verf.) wants to claim that our justification for accepting arithmetical definitions comes from observation of ordinary physical objects. Prominent among such observations are experiences of separating, combining and arranging physical objects. Mill’s idea is to provide a view of arithmetical truth which will exhibit directly the relevance of such observations to our acceptance of the basic laws of arithmetic.} [1980, p. 223]
\end{small}

Entscheidende Gedanken Kitchers dürften auf die Auseinandersetzung mit Mill zurückgehen, insbesondere die Vorstellung, elementare mathematische Tätigkeiten als Idealisierung realer Tätigkeiten aufzufassen:

\begin{small}
\enquote{I shall try to present the thesis that arithmetic provides an idealized description of our physical operations in a way which avoids the unfortunate suggestions of the possible worlds picture.} [p. 230] \smallskip
\end{small}

In [1984] weist Kitcher zunächst jede Form einer platonischen oder
aprioristischen Auffassung von Mathematik zurück und setzt dagegen eine konstruktivistische Sichtweise, die er wie folgt kennzeichnet:

\begin{small}
\enquote{The constructivist position I defend claims that mathematics is an idealized science of operations which we can perform on objects in our environment. Specifically, mathematics offers an idealized description of collecting and ordering which we are able to perform with respect to any object.} [p.12]
\end{small}

\noindent An anderer Stelle heißt es:

\begin{small}
\enquote{Mathematics consists in idealized theories of ways in which we can operate on the world. To produce an idealized theory is to make some stipulations --- but they are stipulations which must be appropriately related to the phenomena one is trying to idealize.} [p. 161]
\end{small} \smallskip

Was sich bei Lakatos abzeichnet, mathematische Theorien als den
naturwissenschaftlichen verwandt aufzufassen, radikalisiert Kitcher, indem er mathematisches Wissen unmittelbar an die Realität bindet. Die Realität wird gewissermaßen zur Quelle dieses Wissens. Er spricht von einer \emph{evolutionary theory} mathematischen Wissens. Ihr Bezug zu Realität kommt in folgenden Sätzen zum Ausdruck:

\begin{small}
\enquote{The principal task of explaining the origins of mathematical knowledge thus becomes one of providing a picture of mathematical reality which will fit with the thesis that our mathematical knowledge can originate in sense perception. [p. 96]

\noindent..... let me face directly the question of how we obtain perceptually based
mathematical knowledge. We observe ourselves and others performing particular operations of collection, correlation, and so forth, and thereby come to know that such operations exist. This provides us with rudimentary knowledge --- protomathematical knowledge, if you like.} [p. 117]
\end{small}

Die letzte Aussage konkretisiert er wie folgt:

\begin{small}
\enquote{Many of the statements of elementary geometry can easily be interpreted by taking them to be part of an idealization which systematizes facts about ordinary physical objects which are accessible to perception.} [p. 124] 
\end{small}

Diese Auffassung der Elementarmathematik betrachtet sie eindeutig als eine empirische Theorie im einleitend genannten Sinne. Die Übereinstimmung von Kitchers Auffassung elementaren mathematischen Wissens zu dem in empirischen Theorien erworbenen wird besonders deutlich bei seinen Ausführungen zur Arithmetik. Dies vermutlich unter dem Einfluß des Millschen Verständnisses. So schreibt er:

\begin{small}
\enquote{We recognize, for example, that if one performs the collective operation called \enquote{making two}, then performs on different objects the collective operation called \enquote{making three}, then performs the collective operation of combining, the total operation is an operation of \enquote{making five}. Knowledge of such properties of such operations is relevant to arithmetic because arithmetic is concerned with collective operations.

As a first approximation, we might think of my proposal as a peculiar form of constructivism. Like the constructivists I hold that arithmetic truths owe their truth ..... to the operations we perform. (...) Unlike most constructivists, I do not think as the relevant operations as private transactions in some inner medium. Instead, I take as paradigms of constructive activity those familiar manipulations of physical objects in which we engage from childhood on. Or, to present my thesis in a way which will bring out its realist character, we might consider arithmetic to be true in virtue not of what \emph{we can do} to the world but rather of what \emph{the world} will let \emph{us do to it}. ..... arithmetic describes those structural features of the world in virtue of which we are able to segregate and recombine objects: the operations of segregation and recombination bring about the manifestation of underlying dispositional traits. [p. 108]

\noindent..... taking arithmetic to be about operations is simply a way of developing the general idea that arithmetic describes the structure of reality. .....
 
\noindent Arithmetic owes its truth not to the actual operations of actual human agents, but to the ideal operations performed by ideal agents. (...) I construe arithmetic as an \emph{idealizing theory}. [p. 109]

\noindent I suggest that we have no way of knowing in advance what powers should be attributed to our ideal subject. Rather the description of that ideal subject and the conditions of her performance must be tested against our actual manipulations of reality.

\noindent That is not to suppose that there is a mysterious being with superhuman powers. Rather, ..... mathematical truths are true in virtue of stipulations which we set down, specifying conditions on the extension of predicates \emph{which actually are satisfied by nothing at all but are approximatey satisfied by operations we perform (including physical operations)}. [p. 110]

\noindent I propose that the view that mathematics describes the structure of reality should be articulated as the claim that mathematics describes the operational activity of an ideal subject. In other words, to say that mathematics is true in virtue of ideal operations is to explicate the thesis that mathematics describes the structure of the world. Obviously, the ideal subject is an idealization of ourselves, but I explicitly reject the epistemological view that we can know a priori the ways in which the idealization should be made. [p. 111]

\noindent At this point it should be clear why, on my account of arithmetic, there is no commitment to the existence of an ideal agent or to ideal operations. Statements of arithmetic, like statements of ideal gas theory, turn out to be vacuously true. They are distinguished from the host of thoroughly uninteresting and pointless vacuously true statements ..... by the fact that the stipulations of actual gases and the fact that the stipulations on the ideal agent abstract from accidental limitations of human agents.} [p. 117, footnote] \smallskip
\end{small}

Kitcher betont Handlungen, die zunehmend allgemeineren Charakter erhalten. Ein Gedanke, den er von Mill übernahm und den aus der Sicht der Psychologie Jean Piaget herausstellt. [1973, pp. 80/81/82]

\begin{small}
\enquote{One central ideal of my proposal is to replace the notions of abstract mathematical objects, notions like that of a \emph{collection} with the notion of a kind of mathematical activity, \emph{collecting}. I have introduced the notion of form, collecting is tied to physical manipulation of objects. One way of collecting all the red objects on a table is to segregate them from the rest of the objects, and to assign them a special place. We learn how to collect by engaging in this type of activity. However, our collecting does not stop there. Later we can collect the objects in thought without moving them about. We become accustomed to collecting objects by running through a list of their names, or by producing predicates which apply to them. Naïvely, we may assume that the production of any predicate serves to collect the objects to which it applies. (...) Thus our collecting becomes highly abstract. We may achieve a hierarchy of collectings by introducing symbols to represent our former collective activity and repeating collective operations by manipulating these symbols. So, for example, corresponding to the set $\{\{a,b\},\{c,d\}\}$ we have a sequence of collective operations: first we collect a and b, then we collect c and d, and, finally we perform a higher level operation on these collections, an operation which is mediated by the use of symbols to record our prior collective activity. As I construe it, the notation \enquote{$\{. . .\}$} obtains its initial significance by representing first – level collecting of objects, and iteration of this notation is itself a form of collective activity.

\noindent Collecting is not the only elementary form of mathematical activity. In addition we must recognize the role of \emph{correlating}. Here again we begin from crude physical paradigms. Initially, correlation is achieved by matching some objects with others, placing them alongside one another, below one another, or whatever. As we become familiar with the activity we no longer need the physical props. We become able to relate objects in thought. Once again, the development of a language for describing our correlational activity itself enables us to perform higher level operations of correlating: notation makes it possible for us not only to talk (e.g.) about functions fom objects to objects (which correspond to certain first – level correlations) but also about functions from functions to functions, and so forth.} [1984, pp. 110/111]
\end{small} 

Die von Kitcher herausgestellte Bedeutung einer geeigneten Sprache, um Zusammenhänge auf der empirischen Ebene angemessen beschreiben zu können, aber auch solche, die sich von ihr abheben (\enquote{introducing symbols}), charakterisiert die Ebene, auf der der \enquote{ideal agent} agiert. Sie ist frei von den Unzulänglichkeiten realer Manipulationen. Kitcher betont, daß die gesamte Mathematik auf dem ersten --- auf empirischer Basis --- erworbenen mathematischen Wissen beruht, das der einzelne erwirbt:

\begin{small}
\enquote{I propose that a very limited amount of our mathematical knowledge can be obtained by observations and manipulations of ordinary things. Upon this small basis we erect the powerful general theories of modern mathematics. [p. 92]

\noindent Thus I see the growth of mathematical knowledge as a process in which an \emph{ur}--practice, a scattered set of beliefs about manipulations of physical objects, gives rise to a succession of multi–faceted practices through rational transitions, leading ultimately to the mathematics of today.} [p. 226]
\end{small} 

Und schon das erste mathematische Wissen ist ohne Einschränkung Mathematik:

\begin{small}
\enquote{Nor is it even necessary to forego the claim, that mathematics studies abstract objects --- \emph{so long as we regard the claim as ultimately interpreted in terms of ideal operations}}. [p. 142], \medskip
\end{small}
 
Daß der
\begin{small}
\enquote{very limited amount of our mathematical knowledge}
\end{small}
 in der Tat in empirischen Theorien erworbenes Wissen meint, belegt folgende Aussage:     

\begin{small}
\enquote{Egyptian surveyors and Babylonian bureaucrats gained mathematical knowledge on the basis of experience. The Greeks showed how their knowledge could be independently grounded, thereby transforming mathematics --- or, perhaps, producing a new science, \emph{pure} mathematics.} [p. 93] \smallskip 
\end{small}

Indem wir den Erwerb erster mathematischer Einsichten an das Verfügen über empirische Theorien binden, beschränken wir den Rahmen, in dem auf diesem Wege mathematische Kenntnisse erworben werden, auf die Elementarmathematik. Mit den Worten von Mac Lane können wir unsern Standpunkt wie folgt formulieren: 

\begin{small}
\enquote{Mathematics starts from a variety of human activities, disentangles from them a number of notions which are generic and not arbitrary, then formalizes these notions and their manifold interrelations. Thus, in the narrow sense, mathematics studies formal structures by deductive methods which, because of the formal character, require a standard of precision and rigor.} [1981, p. 465] \smallskip 
\end{small}

Kitcher beschreibt --- nach unserer Auffassung zutreffend --- die Spezifität des Wissens, das Schüler in empirischen Theorien erwerben. Seine Argumentation und alles weitere, was gesagt wurde, betreffen i.w. ontologische Fragestellungen. Allerdings ist der durch empirische Theorien vermittelte Wissenserwerb nicht so \enquote{selbsterklärend}, wie es die Darstellung Kitchers suggeriert. Dies zeigt ihre strukturalistische Rekonstruktion. Auch elementare empirische Theorien enthalten --- wie alle empirische Theorien --- theoretische Terme/Begriffe. Terme/Begriffe, die in der Theorie, die sie einführt, keine Referenz haben. So ist z.B. in einer empirischen Theorie der Bruchrechnung, die die Addition von Brüchen als das Zusammenfügen von Bruchteilen realisiert, die Multiplikation ein theoretischer Begriff (vgl. [Burscheid und Struve 2009, Kap. III, 1.4]). Die theoretischen Terme/Begriffe stellen sicher, daß die empirischen Theorien nicht nur Wissen enthalten, das aus realem Handeln oder realen Objekten herausgefiltert wird, sondern daß mit ihnen neue Einsichten formuliert werden, die unser Wissen über die WELT erweitern.  

Das Problem der theoretischen Terme/Begriffe in empirischen Theorien ist zwar nicht vermeidbar, ist aber im vorliegenden Fall unabhängig davon, ob die Objekte der Theorie und die Handlungen, die der Lernende mit diesen Objekten vornimmt, von diesem zunehmend als idealisiert betrachtet werden. Gemäß der auf Piaget zurückgehenden Klassifikation befinden sich die Schüler des angesprochenen Alters auf der Stufe der konkreten Operationen, auf der sie bei Unsicherheit hinsichtlich ihrer Schlüsse des Rückgriffs auf Handlungen an realen Objekten bedürfen. Mithin ist die Theorie, die sie erwerben, nach wie vor eine empirische Theorie. \medskip

Wenn das Bild, das Kitcher von Mathematik entwirft, auch stark von den Naturwissenschaften --- insbesondere der Physik (Theorie des idealen Gases) --- geprägt ist, so erhält die Mathematik doch einen eigenständigen Charakter, der frei ist von Einflüssen platonischer Sichtweisen und der empirisch gewonnene Einsichten in eleganter Weise zu integrieren vermag und ihre Bedeutung manifestiert. Ob man wie Kitcher seine Position für das gesamte mathematische Wissen beansprucht, ist nicht zwingend. Beschränkt man sie auf das elementare mathematische Wissen, das im Rahmen der Schule erworben wird, so steht dem nichts entgegen, im Studium zu einem Hilbertschen Verständnis von Mathematik überzugehen. 

Natürlich blieben Kitchers Ausführungen nicht unwidersprochen (vgl. z.B. Thomas Norton–Smith [1992], Sarah Hoffman [2004]), aber die von ihm stark mitgeprägte philosophische Sicht der Mathematik hat zunehmend an Gewicht gewonnen. Es ist nicht zu übersehen, daß die beiden geschilderten Auffassungen --- die einer wachsenden Zahl tätiger Mathematiker und die des Wissenschaftsphilosophen Kitcher --- starke Gemeinsamkeiten haben. Beide haben keinerlei platonische Züge und beide sehen Mathematik als hervorgegangen aus realen Handlungen, stimmen damit überein mit dem von uns in [2009] vertretenen Standpunkt, daß die schulübliche Elementarmathematik im Rahmen empirischer Theorien erworben wird. \smallskip 

Eine interessante Frage ist, wie sich eine bestimmte Auffassung von Mathematik stützen läßt, wie man sie rechtfertigen kann. Daß sie von etlichen Mathematikern und Mathematikphilosophen vertreten wird, legitimiert sie zwar in hohem Maße, ist aber keine Stützung wie sie z.B. eine naturwissenschaftliche Theorie durch ein Experiment erfährt. Da im vorliegenden Falle eine empirische Absicherung entfällt, ist die Frage, ob sich Argumente finden lassen, die nicht nur theoretischen Charakter haben, die man vielleicht als quasi--empirisch bezeichnen könnte.

Das wichtigste Medium, durch das Mathematik weitergegeben wird, ist die Sprache, im Idealfall eine symbolische Sprache aber im Regelfall --- z.B. bei der Vermittlung von Mathematik --- die Umgangssprache. Man kann nun die in der Mathematik verwendeten Begriffe und Schlußweisen der Umgangssprache daraufhin untersuchen, ob und ggf. wie sie an die Realität gebunden ist. Dieser Weg wird seitens einiger Kognitionswissenschaftler beschritten, indem untersucht wird, wie generell abstrakte Begriffe gebildet werden und nach welchen Mustern Denkbezüge erfolgen. Es liegt auf der Hand, daß dabei der Sprache als konkreter Ausprägung unseres Denkens besondere Bedeutung zukommt. \smallskip 

Seit den 1980er Jahren widmen sich u.a. der Philosoph Mark Johnson und der Linguist George Lakoff bevorzugt dieser Aufgabe. In den 1990er Jahren stieß Rafael E. Núñez, ein aus der Psychologie kommender Kognitionswissenschaftler, hinzu. Die Wissenschaftler bezeichnen ihr Arbeitsgebiet als \emph{cognitive linguistics} [Lakoff und Núñez 1997, p. 32], \emph{cognitive semantics} [Johnson 1993], auch als \emph{embodied cognition} [Núñez und Lakoff 1998, p. 87]. \medskip

\noindent \emph{Bemerkung}: Zentrale englischsprachige Begriffe und Kunstwörter werden wir nicht versuchen, ins Deutsche zu übersetzen, da ihr Verständnis dadurch kaum gefördert werden dürfte. \medskip

Das Spektrum --- bevorzugt empirischer --- Forschungen, auf die sich die Wissenschaftler stützen, ist breit gestreut.

\begin{small}
\enquote{Rather, they are based on empirical evidence from a variety of sources, including psycholinguistic experiments (...), generalisations over inference patterns (...), generalisations over conventional and novel language (...), and the study of historical semantic change (...), of language acquisition (...), of spontaneous gestures (...), of American sign language (...), and of coherence in discourse (...).} [Núñez und Lakoff 1998, p. 87] 
\end{small}

Der entscheidende Ansatz ihrer Überlegungen ist die Beobachtung, daß abstrakte Zusammenhänge in sprachlichen Formulierungen als reale Zusammenhänge ausgedrückt werden.

\begin{small}
\enquote{Time, for example, is primarily conceptualised in terms of motion, either the motions of future times toward an observer (\enquote{Christmas is approaching}) or the motion of an observer over a time landscape (\enquote{We’re approaching Christmas})}. [p. 87]
\end{small}

Derartige Metaphern (\enquote{conceptual metaphors}) sind --- so Johnson und Lakoff --- im täglichen Gespräch allgegenwärtig. \\ 
\begin{small}
\enquote{most of our conceptual system is metaphorically structured} [Lakoff und Johnson 1980, p. 106] \\
\end{small}
Die durch die Metaphern vorgenommenen Umdeutungen werden wie selbstverständlich, mühelos und im täglichen Gespräch automatisch verwendet. Sie sind Teil dessen, was sich kognitiv unbewußt vollzieht und beruhen auf struktureller Übereinstimmung mit unserer alltäglichen Erfahrung. 

\begin{small}
\enquote{The system of conventual conceptual metaphor is mostly unconscious, automatic, and used with no noticeable effort, just like our linguistic system and the rest of our conceptual system. \\
Our system of conventual metaphor is \enquote{alive} in the same sense that our system of grammatical and phonological rules is alive; namely, it is constantly in use, automatically, and below the level of consciousness.} [Lakoff 1993, p. 245]
\end{small}    

Metaphern sind nach Auffassung von Johnson und Lakoff aber nicht nur sprachliche Figuren, wie der Begriff gemeinhin verstanden wird, sondern kognitive Werkzeuge:

\begin{small}
\enquote{The metaphor is not just a matter of language, but of thought and reason. The language is secondary. [p. 208] \\
Metaphorical linguistic expressions are surface manifestations of metaphorical thought.} [Lakoff und Núñez 1997, p. 32]
\end{small}

Sie basieren wesentlich auf Erfahrungen im Umgang mit Objekten der Realität, insbesondere dem eigenen Körper. Metaphern bilden die Strukturen der gemachten Erfahrungen -- des \emph{Quellbereichs} (\enquote{source domain}) -- in einen abstrakten Bereich -- den \emph{Zielbereich} (\enquote{target domain}) -- ab, geben so einen Teil der Struktur des Quellbereichs an den Zielbereich weiter.

\begin{small} 
\enquote{Each mapping is a fixed set of ontological correspondences between entities in a source domain and entities in a target domain.} [Lakoff 1993, p. 245]
\end{small}

\noindent In der Formulierung von Johnson:

\begin{small}
\enquote{..... \emph{metaphor} (kursiv durch die Verf.), conceived as a pervasive mode of understanding by which we project patterns from one domain of experience in order to structure another domain of a different kind. So conceived, metaphor is not merely a linguistic mode of expression; rather, it is one of the chief cognitive structures by which we are able to have coherent, ordered experiences that we can reason about and make sense of. Through metaphor, we make use of patterns that obtain in our physical experience to organize our more abstract understanding. Understanding via metaphorical projection from the concrete to the abstract makes use of physical experience in two ways. First, our bodily movements and interactions in various physical domains of experience are structured (...), and that structure can be projected by metaphor onto abstract domains. Second, metaphorical understanding is not merely a matter of arbitrary fanciful projection from anything to anything with no constraints. Concrete bodily experience not only constraints the \enquote{input} to the metaphorical projections but also the nature of the projections themselves, that is, the kinds of mappings that can occur across domains.} [1987, pp. XIV/XV]
\end{small}

Als erstes Beispiel einer Metapher diene die \texttt{The Ideas Are Objects} Metapher, die Aussagen wie \enquote{Ich bekomme ihre zentrale Idee nicht zu fassen}, \enquote{Sie verlor die Schlüsselbegriffe der Theorie aus dem Blick}, \enquote{Halte die Idee einen Augenblick im Gedächtnis}, \enquote{Niemand konnte mit all diesen Ideen gleichzeitig jonglieren} zugrunde liegt. Johnson führt dazu aus:

\begin{small}
\enquote{\texttt{The Ideas Are objects} metaphor is but a tiny part of a massive system of metaphorically defined models we have for the mind, reasoning, understanding, knowledge, truth and language (...). We do not have conceptions of such notions apart from metaphors of these sorts. In other words, (...), our very understanding and cognition of these notions are constituted partly by metaphor. We simply do not and cannot conceptualise \enquote{ideas} apart from using our knowledge of the metaphoric source domain to draw inferences regarding the target domain.} [1993] \smallskip
\end{small}

Daß durch eine Metapher der Zielbereich nur teilweise strukturiert wird, beruht darauf, daß in der Regel nicht alle Komponenten eines Begriffs in den Zielbereich übernommen werden. 

\begin{small}
\enquote{Because concepts are metaphorically structured in a systematic way, e.g., \texttt{Theories Are Buildings}, it is possible for us to use expressions (\emph{construct, foundation}) from one domain (\texttt{Buildings}) to talk about corresponding concepts in the metaphorically defined domain (\texttt{Theories}). What \emph{foundation}, for example, means in the metaphorically defined domain (\texttt{Theory}) will depend on the details of how the metaphorical concept \texttt{Theories Are Buildings} is used to structure the concept \texttt{Theory}. \\
The part of the concept \texttt{Building} that are used to structure the concept \texttt{Theory} are the foundation and the outer shell. The roof, internal rooms, staircases, and hallways are parts of a building not used as part of the concept \texttt{Theory}. Thus the metaphor \texttt{Theories Are Buildings} has a \enquote{used} part (foundation and outer shell) and an \enquote{anused} part (rooms, staircases, etc.). Expressions such as \emph{construct} and \emph{foundation} are instances of the used part of such a metaphorical concept and are part of our ordinary literal language about theories.} [Lakoff und Johnson 1980, p. 52] \medskip 
\end{small} 

Die Strukturen, die durch die Metaphern abgebildet werden, heißen \enquote{image schemata}. Dazu schreibt Núñez: 

\begin{small}
\enquote{Another important finding in cognitive linguistics is that conceptual systems can be ultimately decomposed into primitive spatial relations concepts called \emph{image schemas}. Image schemas are basic dynamic topological and orientation structures that characterize spatial inferences and link language to visual--motor experience (...). As we will see, an extremely important feature of image schemas is that their \emph{inferential structure} is \emph{preserved} under metaphorical mappings. Image schemas can be studied empirically through language (and spontaneous gestures), in particular through the linguistic manifestation of spatial relations.} [2000, p. 11]
\end{small}

\noindent Johnson erläutert dies wie folgt:

\begin{small}
\enquote{..... \emph{image schemata} (kursiv durch die Verf.) are abstract patterns in our experience and understanding that are not propositional in any of the standard sense of that term, and yet they are central to meaning and to the inferences we make.} [1987, p. 2] \smallskip
\end{small}

\noindent \enquote{propositional} sind \enquote{image schemata} im folgenden Sinne: \smallskip

\begin{small}
\enquote{A proposition exists as a continuous, analog pattern of experience or
understanding, with sufficient internal structure to permit inferences.

\noindent I will argue that, because image schemata and their metaphorical extensions are propositional in this special sense, they constitute much of what we call meaning structure and inferential patterns, although they are not finitary. We will see, for instance, that the \texttt{COMPULSION} schema\footnote{\enquote{An actual \texttt{COMPULSION} schema exists as a \emph{continuous, analog} pattern of, or in, a particular experience or cognition that I have of compulsion} [p. 2]}  has internal structure consisting of a force vector (with a certain magnitude and direction), an entity acted upon by the force, and a potential trajectory the entity will traverse. And this structure constraints the way the schema organizes meaning and influences the drawing of inferences in domains of understanding concerned with forces of a certain kind. The main point is that the internal structure of the image schema exists in a continuous, analog fashion within our understanding, which permits it to enter into transformations and other cognitive operations. [pp. 3/4]

..... \emph{image schemata are not rich, concrete images or mental pictures, either}. They are structures that organize our mental representations at a level more general and abstract than that at which we form particular mental images. [pp. 23/24]

Image schemata exist at a level of generality and abstraction that allows them to serve repeatedly as identifying in an indefinitely large number of experiences, perceptions, and image formations for objects or events that are similarly structured in the relevant ways. Their most important feature is that they have a few basic elements or components that are related by definite structures, and yet they have a certain flexibility. As a result of this simple structure, they are a chief means for achieving order in our experience so that we can comprehend and reason about it. [p. 28]

..... image schemata operate at a level of mental organization that falls between abstract propositional structures, on the one site, and particular concrete images, on the other.

The view I am proposing is this: in order for us to have meaningful, connected experiences that we can comprehend and reason about there must be pattern and order to our actions, perceptions, and conceptions. \emph{A schema is the recurrent pattern, shape, and regularity in, or of, these ongoing ordering activities}. These patterns emerge as meaningful structures for us chiefly at the level of our bodily movements through space, our manupilation of objects, and our perceptual interactions. [p. 29]

Schemata are structures \emph{of an activity} by which we organize our experience in ways that we can comprehend. They are a primary means by which we \emph{construct} or \emph{constitute} order and are not mere passive receptacles into which experience is poured. (...) Unlike templates, schemata are flexible in that they can take on any number of specific instantiations in varying contexts. It is somewhat misleading to say that an image schema gets \enquote{filled in} by concrete perceptual details; rather, it must be relatively malleable, so that it can be modified to fit many similar, but different, situations that manifest a recurring underlying structure.} [pp. 29/30] \vspace{3mm}
\end{small}

\noindent Beispiele wichtiger schemata sind

--- das \texttt{Container} schema.
Es besteht aus einem Innern, das ausgezeichnet wird, und als Orientierungspunkt (\enquote{landmark}) dient, einer geschlossenen Grenze und einem Äußeren. Es bildet eine \\
\begin{small}
\enquote{\emph{gestalt} --- a whole that we human beings find more basic than the parts ( ... ); that is, the complex of properties occuring together is more basic to our experience than their separate occurence.} [Lakoff und Johnson 1980, pp. 70/71]
\end{small}

\noindent Das \texttt{Container} schema hat eine Struktur, die vergleichbar ist derjenigen des Satzes vom ausgeschlossenen Dritten. Versteht man Kategorien metaphorisch als Container, so ist jedes Objekt innerhalb bzw. außerhalb des Containers, eine dritte Möglichkeit gibt es nicht. Auch die logische Relation der Transitivität ist eine Projektion der (räumlichen) Transitivität von Containern. 
\smallskip

--- das \texttt{Source--Path--Goal} schema,

\noindent bestehend aus einer Ausgangssituation(= Standort A), einer Ziel--Wunsch--Situation (= Standort B) und einer Handlungsfolge (= Bewegung von A nach B). Es bestimmt alles, was Bewegungen im Raum betrifft, und hat die Struktur einer Vielzahl unterschiedlicher Argumentationsformen, z.B. der \texttt{Purposes Are Destinations} Metapher, dem Verständnis in der Zeit ablaufender Prozesse aber auch (logischer) Begründungen und Schlußfolgerungen, die auf in der Zeit ablaufende
Prozesse projiziert und als solche formuliert werden. [Johnson 1991] \smallskip

\noindent --- Das \enquote{motor--control system} ist derjenige Teil des Nervensystems, der unsere körperlichen Bewegungen leitet. Überraschenderweise ließ sich nachweisen, daß \\
\begin{small}
\enquote{..... \emph{the same neural structure used in the control of complex motor schemas can also be used to reason about events and actions.} ( ... ) neural control systems for bodily motions have the same characteristics needed for rational inference in the domain of \emph{aspect} (kursiv durch die Verf.) --- that is, the structure of events [Lakoff und Núñez 2000, p. 35], the semantic structure of events in general.} [Lakoff und Johnson 1999, p. 41] \\
\end{small}
Eine derartige Struktur wird als \texttt{Aspect} schema bezeichnet. \medskip

Da \enquote{image schemata} (auch) begrifflicher Natur sind, lassen sie sich sog. Transformationen unterwerfen, z.B. in komplexe Zusammenhänge einbinden. 

\begin{small}
\enquote{By \enquote{\emph{transformations}} (kursiv durch die Verf.) I mean such cognitive operations as scanning an image, tracing out the probable trajectory of a force vector, superimposing one schema upon another, and taking a multiplex cluster of entities and contracting it into a homogeneous mass (...).} [Johnson 1987, p. 3] \smallskip
\end{small}

\noindent Zwei Beispiele solcher Transformationen sind \smallskip

--- Conceptual Composition

\noindent Die beiden oben als erste vorgestellten schemata können zum \texttt{Into} schema zusammengefaßt werden, das sich aus dem \texttt{In} schema und dem \texttt{To} schema zusammensetzt. Das \texttt{In} schema ist ein \texttt{Container} schema mit dem Innern als Orientierungspunkt. Das \texttt{To} schema ist ein \texttt{Source--Path--Goal} schema, bei dem das Ziel als Orientierungspunkt dient.

\noindent Es korrespondieren somit (Inneres (1. schema), Ziel (2. schema)) und
folglich (Äußeres (1. schema), Quelle (2. schema)). \smallskip

--- Conceptual Blends

\sloppy \noindent Die begriffliche Verbindung verschiedener kognitiver Strukturen mit festen Übereinstimmungen. Sind die Übereinstimmungen metaphorischer Natur, sprcht man von \emph{metaphorischer Verschmelzung} (\enquote{metaphorical blend}). Ein Beispiel aus der Mathematik ist das Konzept des Zahlenstrahls, das die Descartessche \texttt{Numbers as Points} Metapher mit der \texttt{Line as Set of Points} Metapher zur
\texttt{Numbers are Points on a Line} Metapher verschmelzt. [p. 48] \vspace{4mm}

Lakoff und Núñez übertrugen die oben skizzierten Gedanken auf die Mathematik. Sie stellten sich die Aufgabe, mathematische Ideen aus kognitiver Sicht zu analysieren. Mit \enquote{Ideen} sind hier keine Beweisideen gemeint sondern zentrale mathematische Begriffe wie Zahl, Menge, (algebraische) Operation, wichtige Begriffe der Analysis, das \enquote{Unendliche} o.ä.

\begin{small}
\enquote{It is up to cognitive science and neurosciences to do what mathematics
itself cannot do -- namely, apply the science of mind to human mathematical ideas.} [2000, p. XI] \smallskip
\end{small}

In Übereinstimmung mit entsprechenden Forschungsergebnissen gehen Lakoff und Núñez davon aus, daß bestimmte arithmetische Fähigkeiten --- subitizing (bis zu 4 Elementen), Addition und Subtraktion von 1 zu 1 oder 2 --- schon nach dem ersten halben Lebensjahr beherrscht werden, somit angeborene Fähigkeiten sind. Da der größte Teil unserer Denkvorgänge unbewußt abläuft, folgern Lakoff und Núñez, daß mathematisches Verständnis die gleichen kognitiven Mechanismen verwendet wie das Verstehen in nichtmathematischen Bereichen. Daher gewinnen wir auch mathematische Einsichten mit Hilfe der uns angeborenen kognitiven Fähigkeiten. \\
\begin{small}
\enquote{..... mathematics makes use of our adaptive capacities --- our ability to adapt other cognitive mechanisms for mathematical purposes.} [p. 33] \smallskip
\end{small}

Mathematische Begriffe, die unter Verwendung von Begriffen aus nichtmathematischen Bereichen gebildet werden, sind z.B.

--- der Begriff der Klasse, der sich am Begriff einer Ansammlung von Objekten orientiert,

--- der Begriff der Rekursion, der vom Begriff einer sich wiederholenden Tätigkeit Gebrauch macht,

--- die arithmetischen Operationen im Komplexen, die den Begriff einer Drehung benutzen,

--- der Begriff der Ableitung, der auf die Begriffe Bewegung und Annäherung an eine Grenze zurückgreift. [pp. 28/29]  

Beispiele weiterer mathematischer Begriffe, die durch entsprechende \enquote{image schemata} charakterisiert werden können, sind Zentralisierung, Berührung, Nähe, Ausgewogenheit und Geradheit. \vspace{3mm}

Man macht sich ferner leicht klar, daß sich mit Hilfe der \texttt{Categories Are Containers} Metapher nach Art der Venn--Diagramme die Booleschen Gesetze darstellen lassen. [Lakoff und Núñez 2000, pp. 44/45] \vspace{2mm}

\enquote{Image schemata}, die wesentlich sind für das Verständnis mathematischer Begriffe und die (möglicherweise) universelle, d.h. nicht von der Sprache abhängige “conceptual primitives” bilden, sind

--- das \texttt{Contact} schema (es ist topologischer Natur und kennzeichnet das Fehlen eines Zwischenraums),

--- das \texttt{Support} schema (das die Struktur einer Kraft hat),
            
\noindent (Die Beziehung \enquote{an} (der Wand) beinhaltet diese beiden schemata.)

--- das schon erwähnte \texttt{Container} schema, welches u.a. Basis ist für das Verständnis der Relationen \enquote{innerhalb} und \enquote{außerhalb}. \vspace{2mm}

Für die Entwicklung von Mathematik wichtige Formen von \enquote{conceptual metaphors} sind

--- \emph{grounding metaphors}: Sie gründen unsere mathematischen Ideen auf alltägliche Erfahrung. Bei ihnen ist der Zielbereich mathematisch, der Quellbereich liegt außerhalb der Mathematik. Beispiele sind: Addition als Hinzufügen von Objekten zu einer Kollektion, Mengen als Behälter, Elemente einer Menge als Objekte in einem Behälter.

--- \emph{redefinitional metaphors}: Gewöhnliche Begriffe werden durch mehr technische Versionen ersetzt (So formuliert Georg Cantor \enquote{weniger als} und \enquote{so groß wie} für unendliche Mengen mit Hilfe bijektiver Abbildungen).

--- \emph{linking metaphors}: Metaphern innerhalb der Mathematik. Quell-- und Zielbereich sind mathematisch. Sie erlauben es, einen mathematischen Bereich durch das Begriffssystem eines anderen zu erfassen: Zahlen werden aufgefaßt als Punkte auf einer Linie, geometrische Figuren als algebraische Gleichungen, Operationen auf Klassen als algebraische Operationen.

\enquote{Linking metaphors} sind in mehrfacher Hinsicht vielleicht die interessantesten. Sie sind Teil der Mathematik und sind von zentraler Bedeutung bei der Entwicklung neuer mathematischer Begriffe und Gebiete. Klassische Gebiete wie Analytische Geometrie und Komplexe Analysis verdanken ihnen ihre Existenz. [Núñez 2000, pp. 10/11] \vspace{3mm}

Wir verzichten darauf, weitere kognitive Mechanismen, die Lakoff und Núñez angeben, näher auszuführen. Es dürfte deutlich geworden sein, daß ihre Ausführungen eine Sicht von Mathematik stützen, wie sie i.w. der oben skizzierten Auffassung vieler Mathematiker und des Wissenschaftsphilosophen Kitcher entspricht. Was bei Kitcher der \enquote{ideal agent} leistet, wird von Lakoff und Núñez als Wirkung des \enquote{ embodied mind} nachgewiesen. Allerdings beschränken wir dieses Urteil auf denjenigen Teil der Mathematik, der den Begriff \enquote{unendlich} nicht benötigt. Er umfaßt die heute schulübliche Elementarmathematik.

Daß wir die Elementarmathematik auf diese Weise abgrenzen, hat folgenden Grund. Einer der von Lakoff und Núñez behandelten kognitiven Mechanismen ist \enquote{Metaphors That Introduce Elements}. Er ist entscheidend für die Behandlung des Aktual--Unendlichen in der Mathematik. Die beispielhaft herangezogene \texttt{Love Is a Partnership} Metapher wird so verstanden, als habe sie wenig gemein mit \enquote{Liebe--an--sich}. Der Hinweis, daß Romeo und Julia sich von der Metapher distanzieren würden, da sie neue Vorstellungen in ihre Liebesbeziehung einbringe, darf aber kaum als eine ausreichende Begründung für das gewählte Verständnis angesehen werden.  

Bei der Behandlung des Aktual--Unendlichen ist zunächst zu berücksichtigen, daß Prozesse metaphorisch in der Sprache von Bewegungen gefaßt werden. Prozesse werden gesehen als an--dauernde Bewegungen und \emph{unvollständige} unendliche an--dauernde Bewegungen werden mit der \texttt{Continuous Processes Are Iterative Processes} Metapher gefaßt als sich unendlich oft wiederholende \emph{abgeschlossene} endliche Schritte  
\begin{small}
(\enquote{indefinitely iterating step--by--step processes, in which each step is discrete and minimal}). \\
\enquote{Indefinitely continuous motion is hard to visualize, and for extremly long periods it is impossible to visualize. What we do instead is visualize short motions and then repeat them, thus conceptualizing indefinitely continuous motion as repeated motion.} [Lakoff und Núñez 2000, p. 157]   
\end{small}

Mit Anwendung des Mechanismus \enquote{Metaphors That Introduce Elements} wird dann
jedem sich unendlich oft wiederholenden Prozeß metaphorisch ein \enquote{final resulting state} [p. 159]  zugeordnet. Zur Begründung heißt es:

\begin{small}
\enquote{We \emph{hypothezise} (kursiv durch die Verf.) that the idea of actual infinity in mathematics is metaphorical, that the various instances of actual infinity make use of the ultimate metaphorical \emph{result} of a process without end.} \vspace{3mm}
\end{small}

Die Zuweisung eines eindeutig bestimmten metaphorischen Abschlusses zu jedem unendlichen Prozeß bezeichnen Lakoff und Núñez als Ergebnis der \texttt{Basic Metaphor of Infinity}. Während im Fall der Liebesbeziehung darüber gestritten werden kann, ob die Metapher eine neue Vorstellung einbringt, steht in Fall des unendlichen Prozesses \enquote{Metapher} aus unserer Sicht für eine Setzung, ein Axiom. Lakoff und Núñez geben auch keine Metapher an, die ihre Aussage stützt, die auf die Existenz eines \enquote{final resulting state} hinweist. Dies dürfte auch kaum möglich sein, wenn man so unterschiedliche Prozesse betrachtet wie denjenigen, der die Folge $(\tfrac{1}{n})_{n \in \mathbb{N}}$ hervorbringt, und denjenigen, der eine unendliche (0 -- 1) -- Folge produziert. \\ 
Man fragt sich daher, welchen Stellenwert Analysen mathematischer Begriffe haben, die auf dieser Hypothese fußen. [Lakoff und Núñez 2000, pp. 45/46, 156–161] \smallskip

Trotz dieser Einschränkung ist es ein Verdienst von Lakoff und Núñez, \enquote{image schemata} als Bausteine mathematischer Begriffe identifiziert zu haben, ein zweites, --- nicht minder gewichtiges --- ist, die darauf fußende Auffassung ausgearbeitet zu haben, daß auch mathematisches Denken in weiten Teilen metaphorisch ist.

\begin{small}
\enquote{For the most part, human beings conceptualize abstract concepts in
concrete terms, using ideas and modes of reasoning grounded in the sensory–
motor system. The mechanism by which the abstract is comprehended in terms of the concrete is called conceptual metaphor. Mathematical thought also makes use of conceptual metaphor, as when we conceptualize numbers as points on a line. [p. 5]

Their primary function ist to allow us to reason about relatively abstract domains using the inferential structure of relatively concrete domains.} [p. 42] \medskip
\end{small}

Eine aus unserer Sicht wesentliche Konsequenz ist, daß die Auffassung von Mathematik, die wir einleitend skizziert haben und die offenbar von nicht wenigen Mathematikern geteilt wird, durch kognitionswissenschaftliche Ergebnisse gestützt wird. Diese Auffassung sehen wir daher als gerechtfertigt an. Auch der Standpunkt, daß Schüler ihr mathematisches Wissen im Rahmen empirischer Theorien erwerben, wird damit gestärkt. Zumal es eine mit der obigen Auffassung übereinstimmende mathematik--philosophische Position gibt, in die sich dieses Wissen nahtlos einfügt. \medskip

Was auffällt, ist, daß Lakoff und Núñez nicht die Auswirkungen eines metaphorischen Verständnisses mathematischer Begriffe diskutieren. Daß deren Verständnis besser \enquote{geerdet} wird, dürfte fraglos sein; ein nicht unwesentlicher Punkt.

\begin{small}
\enquote{The discussion of scíentific language suggests that description without explanation is incomplete, and that literal precision without models or metaphors of the underlying mechanisms is sterile. To be useful, the language of science must provide not only precise descriptive information but also cues for how to understand and interprete the information.} [Mayer 1993, p. 566] \vspace{3mm}
\end{small}

Von der Frage des Verständnisses ist aber die Frage der Darstellung mathematischer Inhalte zu trennen. Ob auch Einfluß auf diese genommen wird, ist fraglich. Ein schönes Beispiel liefern Lakoff und Núñez in ihrer Auseinandersetzung mit der Formulierung der Stetigkeit durch Weierstraß. Sie orientieren den Begriff der \emph{natürlichen Stetigkeit} an Euler:\begin{small} \enquote{a curve discribed by freely leading the hand} [Núñez, Edwards und Matos 1999, p. 54].
\end{small}

Nun dürfte sicher sein, daß Euler seine Vorstellung an dem ihm vertrauten Begriff der Kurve orientierte --- wie der Text auch sagt --- und nicht am Begriff der Funktion, dessen vielfältige Ausprägung zu seiner Zeit erst begann. Lakoff und Núñez weisen richtigerweise darauf hin, daß sich Weierstraß’ $\epsilon$--$\delta$ Stetigkeit an der Definition Dedekindscher Schnitte orientiert, daß die \emph{natürliche Stetigkeit} Eulers und die $\epsilon$--$\delta$ Stetigkeit Weierstraß' somit zwei verschieden Begriffe seien --- dies schon deshalb, weil Euler nur die Stetigkeit der \emph{Kurve} im Blick hat, Weierstraß aber die Stetigkeit der die Kurve definierenden \emph{Funktion}. Aber auch für die Graphen reellwertiger Funktionen sind die beiden Begriffe  von unterschiedlicher Stärke. Denn der Graph einer reellwertigen Funktion, der \emph{natürlich stetig} ist (im oben genannten Sinne), ist nicht notwendig gegeben durch eine stetige Funktion im Sinne der Weierstraßschen Definition (Kahle [1971]). \\
In Analogie zu der gebräuchlichen Weierstraßschen Formulierung bringen Lakoff und Núñez den Begriff der \emph{natürlichen Stetigkeit} in Anklang an seine metaphorische Fassung auf folgende Form: \vspace{3mm}

\begin{small}
\noindent \enquote{For instance, lim$_{x\to a}$f(x) = L, could be defined to mean: \\
\noindent \hspace*{4mm}for every $\epsilon$ > 0, there exists a $\delta$ > 0, such that \\
\noindent \hspace*{7mm}as x \emph{moves} toward a and \emph{gets} and \emph{stays} within the distance $\delta$ of a, \\
\noindent \hspace*{5mm} f(x) \emph{moves} toward L and \emph{gets} and \emph{stays} within the distance $\epsilon$ of L}.
\noindent \hspace*{74mm} [Núñez und Lakoff 1998, p. 99] \vspace{1mm}
\end{small}

Diese Formulierung ist deshalb problematisch, weil sie unterstellt, Euler habe die Stetigkeit von Funktionen im Blick gehabt, was aber offensichtlich nicht der Fall war. Selbst wenn man von dieser irrtümlichen Unterstellung absieht, bleibt offen, wie man im konkreten Fall die \emph{natürliche Stetigkeit} einer \emph{Funktion} f an einer Stelle a mit Hilfe dieser Definition feststellt. Ist da nicht die Vorgehensweise der Mathematiker des 17. Jahrhunderts (Leibniz, Bernoulli, ... ) angemessener, die bei den von ihnen untersuchten Kurven auf Existenzbeweise für die beobachteten Phänomene verzichteten, da sie diese \enquote{durch Hinsehen} feststellten, und die Rechnungen nur benutzten, um die interessierenden Stellen
genau zu fixieren?

Somit ist fraglich, ob ein metaphorisches Verständnis mathematischer Begriffe, was --- so kann man die Ausführungen von Lakoff und Núñez lesen --- als \enquote{natürlich} anzusehen ist, den Belangen der Mathematik in allen Punkten gerecht wird. Die Arithmetisierung der Mathematik im 19. Jahrhundert erfolgte zwangsläufig, da andernfalls die Fragen, die der Funktionsbegriff aufwarf, nicht hätten beantwortet werden können. \\
Wie das Beispiel zeigt, ist es zumindest zweifelhaft, ob die \emph{Darstellung} mathematischer Inhalte durch ein metaphorisches Verständnis ihrer Begriffe beeinflußt wird. Dies kommt auch in der folgenden Äußerung zum Ausdruck:

\begin{small}
\enquote{Once up, we throw them (metaphors, die Verf.) away (even hide them) in favor of a formal, logically consistent theory that (with luck) can be stated in mathematical or near-mathematical terms. Their formal models that emerge are shared, carefully guarded against attack, and prescribe ways of life for their users. The metaphors that are added in this achievement are usually forgotten or, if the ascent turns out to be important, are made not a part of science but part of the history of science.} (Bruner [1986, p. 48], zitiert nach Sfard [2008])     
\end{small} \vspace{3mm}

Neben den Fragen nach einem sich wandelnden Verständnis mathematischer Begriffe und nach seinem Einfluß auf die Behandlung mathematischer Inhalte ist natürlich von Interesse, ob und ggf. wie das Erlernen von Mathematik beeinflußt wird. Sieht man von einigen naheliegenden \enquote{Eselsbrücken} wie der \texttt{Categories Are Containers} Metapher ab, so dürfte der Einfluß eher indirekt sein, indem sich das Verständnis von Mathematik wandelt, das sich der Lehrer aneignet. Ob die Erwartungen, die  Núñez e.a. von einem zukünftigen Unterricht hegen, schon allein aus zeitlichen Gründen realistisch sind, ist allerdings die Frage.

\begin{small}
\enquote{..... we schould provide a learning environment in which mathematical ideas are taught and discussed with all their human embodied and social features. Students (and teachers) should know, that mathematical theorems, proofs, and objects are about ideas, and that these ideas are situated and meaningful because they are grounded in our bodily experience as social animals. Providing an understanding of the historical processes through which embodied ideas have emerged can support this aim. This does not mean simply presenting a few names and dates as a prelude to teaching the \enquote{real} mathematics. It means talking about the motivations, zeitgeist, controversies, difficulties, and disputes that motivated and made possible particular in mathematics.} [Núñez, Edwards und Matos 1999, p. 62]  \vspace{3mm}
\end{small}

Einen direkten Einfluß hat die dargestellte kognitionswissenschatliche Konzeption auf die Überlegungen von Anna Sfard. Seit den 1980er Jahren fand die Frage zunehmend Interesse, ob und wie sich ggf. mathematische Begriffe im Verständnis des Lernenden wandeln. Neben Sfard widmete sich insbesondere Ed Dubinsky diesem Problem [1991]. Am Beispiel des Zahl-- und des Funktionsbegriffs zeigte Sfard unter Hinzuziehung historischen Materials, wie schwierig es selbst für wissenschaftlich tätige Mathematiker war, einen Zahl--/Funktionsbegriff zu gewinnen, den man heute als strukturell bezeichnen würde --- ein eigenständiges Element einer komplexen Struktur, obwohl der Zahlbegriff (operational) und der Funktionsbegriff (intuitiv) in der Praxis wie selbstverständlich gehandhabt wurden [1991, 1992].

Der Zahlbegriff ist ein Beispiel der großen Klasse mathematischer Begriffe, die zunächst operational auftreten. Um den Übergang zu einem strukturellen Verständnis dieser Begriffe in ihrem Sinne darzustellen müssen wir auf Sfards Sichtweise von Bedeutungsgewinnung eingehen.  

Im Zentrum ihrer Überlegungen stehen die Ausführungen von Lakoff und Núñez. In enger Anlehnung an diese Autoren führt sie hinsichtlich des \enquote{mechanism of metaphorical construction} aus:

\begin{small}
\enquote{..... the vehicle which carries our experimentally constructed knowledge is an \emph{embodied schema} (known also as an \emph{image schema}). Johnson defines embodied schemata as \enquote{structures \emph{of an activity}} by which we organize our experience in ways that we can comprehend. They are a primary means by which we \emph{construct} or \emph{constitute} order and are not mere passive receptacles into which experience is poured. ..... an embodied schema is what epitomizes, organizes, and preserves \enquote{for future use} the essence of our experience and, as such, it is our tool for handling the multifarious physical and intellectual stimuli with which are faced throughout our lives.   

..... embodied schemata (...) are image--like and embodied, embodied in the sense that they should be viewed as analog reflections of bodily experience rather than as factual statements we may wish to check for validity. The non--propositional nature of embodied schemata makes it difficult, sometimes impossible, to describe them in words. [1994, p. 46]

If embodied schemata cannot be viewed as the mental counterpart of a system of factual statements, the question arises about the cognitive means by which such schemata are handled. Here again, misled by our previous knowledge, we may easily slip into an oversimplified, distorted version. Mental images seem to be the natural alternative to the propositional structure. The idea that an embodied schema is, in fact, a mental image is even more convincing in view of the fact that both these cognitive structures have the same leading characteristics: they are \emph{analog and holistic}. ..... whereas a mental image is always an image of something concrete and and is therefore full of details ( ..... ), an embodied schema is general and malleable. It is but a skeleton with many variable parts which, being undetermined, cannot be visualized.} [p. 47] \vspace{3mm}
\end{small} 

Um die Überlegungen von Sfard zu konkretisieren werden wir sie in ein ausführlich dargestelltes Beispiel einbeziehen. Wir wählen eine Einführung der reellen Zahlen --- damit das Problem der irrationalen Zahlen. Dieses Thema ist deshalb für von besonderem Interesse, weil es einen entscheidenden Wechsel im Zahlverständnis beinhaltet. Denn aus unserer Sicht ist das so gut wie stets auftretende Verständnisproblem des Schülers darin begründet, daß die bis dato erworbenen Zahlbegriffe (natürliche, ganze, rationale Zahlen) in seinem Verständnis \emph{empirische} Begriffe sind, was heißen soll, daß die Zahlen und die Operationen mit ihnen in empirischen Modellen realisiert werden können. Sich von einem empirischen Zahlverständnis zu lösen dürfte zu den Hauptschwierigkeiten eines Studienanfängers gehören, der eine Analysisvorlesung besucht, die von ihm mehr verlangt als das Beherrschen gewisser Techniken. Denn mit einem empirischen Zahlverständnis sind tragende Begriffe der Analysis (z.B. alle solchen, die auf Konvergenz fußen) nicht zu verstehen. Es geht uns also nicht darum zu betonen, daß der Begriff der irrationalen Zahl schwer verständlich ist --- dies wäre ein Gemeinplatz, sondern wir heben darauf ab, daß die entscheidende Schwierigkeit ist, den vertrauten empirischen Zahlbegriff durch einen --- im noch zu bestimmenden Sinne --- theoretischen Begriff zu ergänzen. Das Vorgehen in vielen Anfängervorlesungen, die reellen Zahlen axiomatisch einzuführen, löst dieses Problem nicht, da den Hörern die Sinnhaftigkeit dieses Vorgehens nicht klar sein dürfte, und es auch nicht ihr Verständnisproblem aufgreift, wie wir zeigen werden. \medskip

Eine der praktikabelsten Einführung der reellen Zahlen, die in den letzten Jahrzehnten vorgelegt wurden, ist der Entwurf von Heinrich Bürger und Fritz Schweiger [1973]. Er orientiert sich am Axiom der Existenz des Supremums und macht sich zunutze, daß der Schüler mit der linearen Ordnung der rationalen Zahlen vertraut ist. An diesem Beispiel, das dem Verständnis von Schülern besonders entgegenkommt, wollen wir zeigen, wie ausgeprägt sich für sie der Begriff der irrationalen Zahl von den bislang erworbenen empirischen Zahlbegriffen unterscheidet. Formal stellen wir den Entwurf nur soweit dar, wie es die hier verfolgte Intention erfordert. Wie die Autoren sprechen wir im folgenden als Adressaten den Schüler an. \smallskip 

\sloppy Ausgehend von der Struktur ($\mathbb{Q}^+;\leqq;+,\cdot$) stützen Bürger und Schweiger die Konstruktion auf bestimmte Teilmengen von $\mathbb{Q}^+$, die sie als rationale Strecken bezeichnen. Eine \emph{rationale Strecke} S definieren sie wie folgt:\\

\begin{tabular}{l l c}
$S\subseteq \mathbb{Q}^+ \: $ mit &(i)& $\: S  \ne \emptyset \wedge S \ne \mathbb{Q}^+$ \medskip \\
& (ii) & $ x \in S \wedge y \in \mathbb{Q}^+ \wedge y < x \rightarrow y \in S$ \medskip \\
& (iii) & $ x \in S \rightarrow \underset{z \in S}{\bigvee}(x < z)$ \bigskip \\
\end{tabular}

\noindent Strecken im Sinne von Bürger und Schweiger sind demnach offene Intervalle in $\mathbb{Q}^+$ mit dem linken Endpunkt 0.

Wie stellen sich diese Objekte aus der Sicht des Schülers dar? Da Schüler in der Regel nur $\pi$ und $\sqrt{p}$ (p prim) --- vielfach nur $\pi$ und $\sqrt{2}$ --- als nichtrationale Zahlen kennen, beschreibt die Definition \emph{für sie} Objekte, die sie als bei 0 beginnende \emph{geometrische} Strecken auffassen, d.h. \emph{im Verständnis der Schüler} handelt die Konstruktion von geometrischen Objekten --- empirisch realisierbaren Elementen. Wenn die Strecken Realisanten haben, die für Schüler empirisch verifizierbar sind, sind für sie diese Verifikationen nicht nur geeignete Veranschaulichungen, sondern der gesamte Konstruktionsprozeß erhält für sie empirischen Charakter (vgl. Struve [1990]), ist Teil einer \emph{empirischen} Theorie und nicht einer mathematischen im heutigen Verständnis. Der Konstruktionsprozeß dürfte bei einer unterrichtlichen Behandlung folglich auf zwei völlig verschiedenen Ebenen ablaufen: \\
--- Der Ebene des Lehrers, auf der dieser einen begrifflich präzisen mathematischen Theorieteil formuliert,\\
--- der Ebene des Schülers, auf der die Objekte, von denen die Konstruktion handelt, empirischer Natur sind, auf der der Schüler  eine Theorie entwickelt, die den Theorien der experimentellen Naturwissenschaften gleicht, eine \emph{empirische Theorie}. \medskip

\noindent \emph{Bemerkung}: Wenn wir sagen, der Schüler entwickelt eine empirische Theorie, so wollen wir damit ausdrücken, daß der Schüler die Objekte und die Relationen zwischen ihnen so handhabt, als bestimme eine solche Theorie sein Handeln, nicht, daß er formal über ein solches Konstrukt verfüge. Daß eine solche Sichtweise gerechtfertigt ist, zeigen die eingangs genannten  Untersuchungen von Psychologen. \medskip

\sloppy Wenn die Konstruktion für Schüler auf der empirischen Ebene abläuft, ist es nur konsequent, sie auch als eine empirische Theorie zu beschreiben. Eine solche Theorie wird im folgenden angegeben. \smallskip 

Im Gegensatz zur Darstellung einer mathematischen Theorie, in der alle Begriffe den gleichen Status haben --- es sind Variable, und deren Darstellung, sofern die Setzungen (Axiome) der Theorie bekannt sind, insofern unproblematisch ist, als man nur die logische Abfolge der Aussagen zu berücksichtigen hat, ist die Frage offen, wie sich eine empirische Theorie \emph{darstellen} läßt, wie sie geeignet formuliert werden kann. Denn sie enthält neben Beg
riffen, die ihre Referenz innerhalb der Theorie haben, auch solche, für die dies nicht gilt --- die theoretischen Terme/Begriffe der Theorie (s.o.). Sie sind es, die dem Lernenden besondere Verständnisschwierigkeiten bereiten. \smallskip

Zur Darstellung einer empirischen Theorie greifen wir eine Theorieform auf, die die Strukturalisten entwickelt haben (Stegmüller [1973, 1986], Balzer [1982]). Eine kurze Einführung findet die Leserin/der Leser bei Burscheid und Struve [2009].

In dieser Theorieform sind Theorien keine Satzklassen (\emph{statement view}) sondern werden zerlegt in eine mathematische Grundstruktur und eine Menge \emph{intendierter Anwendungen}, solcher Probleme, zu deren Lösung die Theorie herangezogen werden soll (\emph{non--statement view}). Die mathematische Grundstruktur ist gegliedert und läßt in einer systematischen Weise die Konstruktion der Theorie erkennen. Man axiomatisiert die Theorie durch Angabe eines mengentheoretischen Prädikates. Aus Gründen besserer Praktikabilität benutzt man keine formale Sprache sondern die auf Patrick Suppes zurückgehende \enquote{informelle mengentheoretische Axiomatisierung} (Stegmüller [1985, S. 39]). Diese Vorgehensweise ist aus der Mathematik geläufig. Als Beispiel sei die folgende Definition genannt: Ein Tripel (P, L, I) heißt heißt \emph{affine Ebene}, wenn P und L nicht--leere Mengen sind, I $\subseteq P \times$ L eine Relation ist zwischen Elementen aus P und Elementen aus L, und wenn folgende Axiome gelten ... . Die Tripel (P, L, I), die die Axiome erfüllen, werden \emph{Modelle affiner Ebenen} genannt.

\sloppy Ganz analog definiert man das Prädikat \enquote{ist eine \textsf{ Klassische Partikelmechanik}}: Ein Tupel (P, T, S, m, f) heißt \emph{\textsf{Klassische  Partikelmechanik}}, wenn P eine endliche, nicht--leere Menge ist (die Menge der \enquote{Partikel}), T ein Intervall reeller Zahlen (ein Intervall von \enquote{Zeitpunkten}), s eine Funktion von P $\times$ T nach $\mathbb{R}$ (die Massenfunktion) und f eine Funktion von $P \times T \times \mathbb{N}$ nach $\mathbb{R}^3$ (die Kraftfunktion), so daß das 2. Newtonsche Gesetz gilt.

Die Tupel, die eine durch die Axiome charakterisierbare mathematische Struktur haben (in dem genannten Beispiel die \textsf{Klassischen Partikelmechaniken}) nennt man \emph{Modelle} der Theorie. Die Strukturen, von denen es sinnvoll ist zu fragen, ob sie Modelle sind, werden \emph{potentielle Modelle} genannt. Die Modelle sind genau diejenigen potentiellen Modelle, die die Axiome der Theorie erfüllen. Die potentiellen Modelle der \textsf{Klassischen Partikelmechanik} sind die Tupel (P, T, s, m, f), deren Komponenten wie oben definiert sind, für die aber nicht notwendig das 2. Newtonsche Gesetz gilt.

Neben den potentiellen Modellen werden noch partiell potentielle Modelle (kurz: \emph{partielle} Modelle) eingeführt. Um den Unterschied zwischen potentiellen und partiellen Modellen zu verdeutlichen, muß man auf die Sprache der Theorie eingehen. Interessante Theorien erweitern die Sprache durch neue Funktions-- und Relationszeichen, die für die Theorie in dem Sinne charakteristisch sind, daß sie erst durch diese eine Interpretation erfahren. Der Teil der Sprache, der ohne die Theorie T interpretierbar ist, heißt (T--)\emph{vortheoretisch}, der Rest heißt (T--)\emph{theoretisch}. Die möglichen Anwendungssituationen von T gehören zu den partiellen Modelle von T. Sie sind Strukturen, die (T--)vortheoretisch beschrieben werden können, bevor die Theorie aufgestellt ist. Die die Modelle charakterisierenden Axiome enthalten mindestens einen T--theoretischen Begriff. Gehörten alle Begriffe der vortheoretischen Sprache an, würde die zu formulierende Theorie T keine Aussage machen, die nicht ohne die Theorie einzusehen wäre. (In einer mathematischen Struktur gibt es natürlich keine theoretischen Begriffe in diesem Sinne.) Da die Modelle in der erweiterten Sprache definiert werden, ist es nicht möglich zu fragen, ob ein partielles Modell ein Modell der Theorie ist.

Die partiellen Modelle der \textsf{Klassischen Partikelmechanik} sind die Tupel (P, T, s), wobei die Komponenten wie oben definiert sind. 

Bezeichnen wir mit M bzw. M$_p$ bzw. M$_{pp}$ die Menge aller Modelle bzw. aller potentiellen Modelle bzw. aller partiellen Modelle einer empirischen Theorie T, so ist die Struktur (M, M$_p$, M$_{pp}$) noch ohne einen Bezug zur WELT, über die eine empirische Theorie ja Aussagen machen möchte. Man betrachtet daher eine Menge I \emph{intendierter Anwendungen} von T, d.h. solcher Problembereiche, die durch T erklärt werden sollen. Die Menge I wird nicht extensional sondern durch \emph{paradigmatische Beispiele} festgelegt. Für die \textsf{Klassische Partikelmechanik} sind dies z.B. das Sonnensystem oder die Pendelbewegungen.

Paradigmatische Beispiele, die I beschreiben, grenzen das Umfeld ab, in dem die zu entwickelnde Theorie angewendet werden soll. Wählt man eine Problemsituation außerhalb dieses Umfeldes, so ist diese der Theorie nicht zugänglich. Ein schönes Beispiel liefert die durch Perlen-- oder Plättchenmengen realisierte Addition (kleiner) natürlicher Zahlen. In den paradigmatischen Beispielen, die man den Kindern vorführt, sind die zu vereinigenden Mengen disjunkt. Solange diese Voraussetzung aber von den Kindern nicht als notwendig eingesehen wird, vereinigen sie auch nicht--disjunkte Mengen, was von der Theorie aber nicht abgedeckt wird.

Wie stellt man nun den Bezug her zwischen der Menge M der Modelle der Theorie T  und einer Menge I intendierter Anwendungen? Traditionellerweise würde man diese \emph{empirische Behauptung der Theorie} folgendermaßen formulieren: Die intendierten Anwendungen von T sind Modelle von T (d.h. I $\subseteq$ M). Um dem \enquote{Problem der theoretischen Terme/Begriffe} gerecht zu werden, muß man zum \enquote{Ramsey -- Substitut} dieser Behauptung übergehen. Letzteres kann man wie folgt formulieren: \emph{Die intendierten Anwendungen von T sind partielle Modelle, die man durch Hinzufügen geeigneter theoretischer Größen/Begriffe zu Modellen von T ergänzen kann}. Die Überprüfung dieser Aussage zerfällt in zwei Teile: Man muß erstens zeigen, daß I in M$_{pp}$ enthalten ist und zweitens, daß die Elemente von I ($\subseteq$ M$_{pp}$) sich zu Modellen ergänzen lassen. Der zweite Schritt ist in dem Sinne unproblematisch, daß man weiß, was man zu tun hat. Im ersten Schritt steckt, wie Wolfgang Balzer formuliert, \enquote{fast die gesamte Problematik der Erkenntnistheorie} [1982, S. 288]: Wie kommt man von realen Systemen, die man in der Welt aufweisen kann, zu theoretischen Strukturen, in unserem Fall zu partiellen Modellen? Diese Frage wird durch den vorgestellten Formalismus natürlich nicht beantwortet. \\ 

\noindent \emph{Bemerkung}: Hier wird deutlich, daß das strukturalistische Theorienkonzept lediglich der \emph{Darstellung} empirischer Theorien dient. Es macht keine philosophischen Aussagen. Der ontologische Status der Begriffe wird durch die Darstellung nicht berührt. Die Unterscheidung der Begriffe einer empirischen Theorie in vortheoretische und theoretische orientiert sich ausschließlich daran, wie der jeweilige Begriff in diese Theorie einbezogen wird.  So ist kein Begriff schlechthin theoretisch. Derselbe Begriff kann hinsichtlich der einen Theorie theoretisch, hinsichtlich einer anderen aber vortheoretisch sein, z.B. wenn er einer schon bekannten Theorie angehört. \\ 

Die bisherige Beschreibung empirischer Theorien ist noch in einem wesentlichen Punkt zu ergänzen. Empirische Theorien haben i.a. nicht eine einzige sozusagen \enquote{kosmische} Anwendung, sondern zahlreiche verschiedene Anwendungen. Diese können sich überschneiden, wie etwa in der Newtonschen Theorie die Systeme Sonne --- Erde und Erde --- Mond. In diesen beiden Anwendungen wird man zweckmäßigerweise der Erde die gleiche Masse zuordnen. Ein logischer Zwang dazu besteht aber nicht, denn es handelt sich ja um voneinander verschiedene Anwendungen. Solche \emph{Querverbindungen} zwischen verschiedenen Anwendungen definiert man extensional als Mengen potentieller Modelle, inhaltlich gesprochen als die Mengen derjenigen potentiellen Modelle, die zu je zweien die jeweils betrachtete Bedingung erfüllen. Bezeichnet C die Menge aller Querverbindungen von T, so können wir definieren: \smallskip \\
Eine \emph{empirische Theorie} ist ein Quintupel (M$_{pp}$, M$_p$, M, C, I) mit
\par
\begingroup
\leftskip=2cm
\noindent
M$_{pp}$: Menge aller partiellen Modelle \\
M$_p$   : Menge aller potentiellen Modelle \\
M$\:$ \,: Menge aller Modelle \\
C $\;$  : Menge aller Querverbindungen \\
I \quad : eine Menge intendierter Anwendungen 
\par
\endgroup \smallskip
\noindent Die oben formulierte empirische Behauptung einer Theorie muß man durch den folgenden Nachsatz ergänzen: \enquote{ ... \emph{und zwar so, daß die Querverbindungen erfüllt sind}.} 

Soweit eine knappe Skizze der theoretischen Struktur, deren wir uns bedienen. Wir betonen nochmals, daß die strukturalistische Theorieform, lediglich der \emph{Darstellung} empirischer Theorien dient, nicht aber deren Charakter berührt, insbesondere nicht den ontologischen Status ihrer Begriffe beeinflußt. \smallskip

Das Gesagte läßt sich wie folgt zusammenfassen: Empirische Theorien haben ontologische Bindungen. Diese kommen darin zum Ausdruck, daß die intendierten Anwendungen zu den partiellen Modellen gehören und zu Modellen ergänzt werden können. Man kann dies auch so ausdrücken, daß über die intendierten Anwendungen der Kontext der Theorie zu ihrem Bestandteil wird. (Dies zeigt erneut, daß es zutreffend ist, Schülertheorien als empirische Theorien aufzufassen.) Die Theoriebildung in empirischen Theorien erfolgt so, daß man versucht, verschiedene Anwendungen, deren Erklärung man intendiert, durch die gleiche Gesetzmäßigkeit zu beschreiben. Dazu bedarf man der theoretischen Begriffe. Diese sind solche, in denen nicht nur in der WELT aufweisbares Wissen ausgedrückt wird. \smallskip 

Wir kommen nun zur Definition der Theorie. Wir bezeichnen sie als Theorie--Element T$_{\mathbb{R^+}}$. T$_{\mathbb{R^+}}$ beschreibt also den Konstruktionsprozeß \emph{aus der Sicht des Schülers}. \smallskip

Die Angabe des Theorie--Elementes T$_{\mathbb{R}^+}$: Als erstes \smallskip \\
die \emph{partiellen Modelle}: \smallskip

M$_{\text{pp}}$(T$_{\mathbb{R}^+})$ = $\langle \mathbb{S} \rangle$ mit \smallskip 

$\mathbb{S}$: die Menge aller geometrischen Strecken $\bar{S}, \bar{T}, ... ,$ angetragen bei 0 in Richtung einer positiven Zahlengeraden. Die Endpunkte einer jeden Strecke werden als ihr nicht zugehörig betrachtet.\\

Die an die Strecken gestellte Bedingung, ihre Endpunkte nicht als ihnen zugehörig zu betrachten --- was dem empirischen Charakter des Theorie--Elementes T$_{\mathbb{R}^+}$ zu widersprechen scheint --- läßt sich am Messen einer Strecke erläutern. Dazu fassen wir im folgenden --- wie Bürger und Schweiger --- die Strecken als Punktmengen auf und identifizieren die Punkte mit ihren Koordinaten auf einer Zahlengeraden, sodaß wir die Strecken auch als Intervalle auffassen können. Betrachtet man die Strecke (0,a) mit a > 0 im Sinne der Definition und die \enquote{Strecke} [0,a], so erhält man beim Messen ihrer Längen das gleiche Ergebnis, d.h. auf der empirischen Ebene unterscheiden sich (0,a) und [0,a] nicht. Die Endpunkte der Strecken außer Betracht zu lassen, hat also keine empirische Relevanz und damit keinen Einfluß auf den Konstruktionsprozeß des Schülers. \\

Bürger und Schweiger behandeln als erstes das Rechnen mit rationalen Strecken --- was nach dem Gesagten in den partiellen Modellen von T$_{\mathbb{R}^+}$ möglich ist, wodurch der empirische Charakter der Konstruktion für den Schüler noch verstärkt werden dürfte. Denn Intervalle der Form (0,a) mit a $\notin \mathbb{Q}^+$ dürften ihm kaum begegnet sein. Wir beschränken uns hier auf die Behandlung der Ordnungsrelation für rationale Strecken, da sich an ihr das Hauptproblem der Schüler verdeutlichen läßt. \smallskip 

\noindent Die \emph{potentiellen Modelle}: \medskip

\begin{tabular} {l l l}
M$_{\text{p}} \big($T$_{\mathbb{R}^+} \big )$ & = & $\big\langle \langle \mathbb{S}, \big \lbrace \underset{\bar{S} \in \mathbb{S'}}{\bigcup} \bar{S} \: \arrowvert \:
\mathbb{S}' \subseteq \mathbb{S} \big \rbrace \rangle \big \rangle$ mit \\

& (i) & $\mathbb{S} \in $ M$_{\text{pp}}\big($T$_{\mathbb{R}^+} \big )$ \medskip \\

& (ii) & $\,\emptyset \neq \mathbb{S}' \,  \subseteq \mathbb{S}: \underset{\bar{S} \in \mathbb{S}'}{\bigcup}\bar{S} \underset{\text{def}}{=}  
\big \lbrace x \, \arrowvert \underset{\bar{S} \in \mathbb{S}'}{\bigvee} (x \in \bar{S}) \big \rbrace $ \bigskip \\
\end{tabular}

\sloppy (ii) erweitert die mengentheoretische Sprache des Schülers. Im Theorie--Element T$_{\mathbb{R}^+}$ sind Vereinigungen über unendlich viele Mengen erforderlich. Da unendlich viele Operationen auf der empirischen Ebene der partiellen Modelle nicht realisierbar sind, nehmen wir $\underset{\bar{S} \in  \mathbb{S}'}{\bigcup}\bar{S} \; (\mathbb{S}' \, \subseteq \mathbb{S})$ als theoretische Terme in die Sprache von T$_{\mathbb{R}^+}$ auf. Sie sind zunächst nur mathematische Zeichen, die keinen Objektcharakter besitzen, d.h. sie lassen sich nicht als Elemente einer komplexen Struktur auffassen. \vspace{3mm}    

Es ist dies der Punkt, an die Überlegungen von Sfard anzuknüpfen. Den Übergang eines zunächst operational verstandenen Begriffs oder auch eines operational interpretierbaren Zeichens --- wie des Zahlbegriffs oder des Zeichens $\underset{\bar{S} \in \mathbb{S'}}{\bigcup} \bar{S}$ --- zu einem strukturellen Verständnis betrachtet sie als einen Wechsel seines
\enquote{embodied schema}. Sie gliedert diesen Wechsel in drei Stufen. Auf der ersten, der \emph{interiorization}, wird der Lernende auf empirischer Ebene --- also in den partiellen Modellen --- mit dem Prozeß vertraut, der dem zu bildenden Begriff zugrundeliegt. Beim vorliegenden Beispiel bedeutet dies im endlich Fall, durch Vergleichen eine jeweils längere Strecke zu bestimmen. So gewinnt der Lernende die Einsicht, daß jede rationale Strecke \enquote{nach oben} fortgesetzt werden kann. Dieses sukzessive Vergleichen und Auswählen der längeren Strecke ist der Prozeß, den der Schüler als mentales Konstrukt, als \enquote{operational schema}, ausbildet: Er erhält ein operationales Verständnis des Begriffs/Zeichens.

Die zweite Stufe, \emph{condensation}, betrachtet den Prozeß des fortwährenden Verlängerns als ganzes, zerlegt ihn nicht mehr in einzelne Bestandteile. 

Die dritte Stufe, \emph{reification}, charakterisiert Sfard als den Übergang von einem \enquote{operational} zu einem \enquote{structural embodied schema}. Dazu führt sie aus: 

\begin{small}
\enquote{An \emph{operational schema} brings into the domain of abstraction a metaphor of doing, of operating on certain objects to obtain certain other objects. As such, it is a schema of action.

The \emph{structural embodied schema}, on the other hand, conveys a completely different ontological message --- a message about a permanent, object--like construct which may be acted to produce other constructs. The advantage of the latter type of schema over the former is that it is more integrative, more economical, and manipulable, more amenable to holistic treatment (...).} [1994,p. 53] \medskip
\end{small}

Zieht man das strukturalistische Begriffssystem zur Darstellung empirischer Theorien heran, so kann man die erste Stufe (interiorization) auffassen als eine Kennzeichnung der Ebene des Wissens, die die partiellen Modelle einer empirischen Theorie repräsentieren. Auf dieser Ebene kann man das Wissen als betont handlungsorientiert betrachten 
\begin{small}
(\enquote{operating on certain objects to obtain other objects}),  
\end{small}
während dieses Operieren auf der Ebene der potentiellen Modelle als eigenständiger Prozeß (condensation) anzusehen ist, in die potentiellen Modelle aufgenommen als undefinierter Term/Begriff. In den Modellen erhalten diese Terme/Begriffe durch eine \enquote{geeignete} Operationalisierung --- \enquote{geeignete} Verknüpfung mit den vortheoretischen Begriffen --- ihre Bedeutung (verstanden im Sinne von Wittgensteins \enquote{Bedeutung (eines Terms/Begriffs) $\widehat{=}$ Gebrauch (der/des den Term/Begriff enthaltenden Formel/Gesetzes}) in der Theorie (reification). Wissen über Prozesse und Handlungen bekommt die Form von Aussagen über Zustände und Objekte. \smallskip 

Im hier behandelten Beispiel bedeutet reification --- der Wechsel des
embodied schema, daß die $\underset{\bar{S} \in \mathbb{S}'}{\bigcup} \bar{S}$ nach M$_{\text{p}} \big($T$_{\mathbb{R}^+} \big )$(ii) Objektcharakter erhalten. \smallskip

Von einem Wechsel des \enquote{embodied schema} eines Begriffs ist bei Johnson die Rede, bei Lakoff und Núñez von unterschiedlichen Weisen, einen Begriff zu fassen. So ließe sich im vorliegenden Beispiel der Wechsel im Begriffsverständnis aus der Sicht von Lakoff und Núñez auch wie folgt interpretieren:

Die erste Stufe wird aufgefaßt wie bei Sfard. Dann wird davon Gebrauch gemacht, daß wir Prozesse sowohl dynamisch (verlaufend in der Zeit) wie auch statisch (wie Behälter, Bahnen von Bewegungen, physikalische Objekte) auffassen.

\begin{small}
\enquote{Processes, as we ordinarily think of them, extend over time. But in mathematics, processes can be conceptualized as atemporal. For example, consider Fibonnacci sequences, in which the n + $2^{\text{nd}}$ term is the sum of the n$^{\text{th}}$ term and the n + 1$^{\text{th}}$ term. The sequence can be conceptualized either as an ongoing infinite process of producing ever more terms or as a thing, an infinite sequence that is atemporal. This dual conceptualization, as we have seen, is not special to mathematics but part of everyday cognition.

..... the difference will not matter, since we all have conceptual mechanisms (...) for going between static and dynamic conceptualizations of processes.} [2000, p. 163] \medskip
\end{small} 

\noindent \emph{Bemerkung}: Eddie M. Gray und David O. Tall haben das Ergebnis von Sfard aufgegriffen und in geschickter Weise erweitert. Sie nutzen aus, daß der Prozeß und der das Prozeßergebnis bezeichnende Begriff im arithmetisch--algebraischen Bereich die gleiche symbolische Darstellung haben können. 4 + 5 bezeichnet eine Addition wie auch eine Summe, eine Darstellung der Zahl 9. Gray und Tall führen den Begriff \emph{elementary procept} ein als Verschmelzung dreier Komponenten, eines \emph{Prozesses}, der ein mathematisches \emph{Objekt} hervorbringt und eines \emph{Symbols} als Darstellung für beide. Hat man mehrere elementary procepts mit demselben Objekt, so sprechen sie von einem \emph{procept}. Die Verbindung von begrifflichem und prozeduralem Denken beim Umgang mit procepten nennen sie \emph{procepturales Denken}. [1994]

Der Begriff \emph{elementary procept} erfaßt genau den Wandel, den ein operational definierter Begriff eines partiellen Modells durchlaufen muß, um in der Symbolsprache der Mathematik in ein \enquote{zugehöriges} potentielles Modell aufgenommen zu werden. Er konkretisiert gewissermaßen die Aussagen von Sfard. \medskip

Für diejenigen Begriffe einer empirischen Theorie, die in der Theorie eine Referenz haben, wurde soweit gezeigt, welchen Bezug sie zu den partiellen und den potentiellen Modellen haben. Noch offen ist, wie die theoretischen Elemente/Begriffe in die informelle Semantik der Theorie einbezogen werden. In der Formulierung der potentiellen Modelle sind sie nur mathematische Zeichen. Die \emph{informelle Semantik} (Stegmüller [1979, S. 479]) (Semantik verstanden im Sinne von Wittgensteins \enquote{Bedeutung $\widehat{=}$ Gebrauch}), die sich u.a. in der Stufung der Modellformen ausdrückt, läßt den Aufbau und in einem systematischen Sinne die Entwicklung der Theorie erkennen. Natürlich wird nicht unterstellt, daß der Schüler eine empirische Theorie gemäß dieser Systematik erwirbt, aber sie legt alle Schritte offen, die zum Erwerb der Theorie erforderlich sind. Defizite des Schülers können vor dem Hintergrund der Systematik leichter identifiziert werden. Auch die Unterscheidung der Begriffe spielt hier hinein. Vortheoretische Begriffe --- sofern sie nicht schon mit einer anderen Theorie erworben wurden --- lassen sich operational oder ostensiv erwerben. Für theoretische Begriffe gilt dies nicht. Sie haben in der Theorie keine Referenzen, kommen daher nicht in den partiellen Modellen vor. Sie werden --- wie schon gesagt --- als undefinierte Terme in den potentiellen Modellen eingeführt 

Gemäß der Absicht, den Umgang mit der WELT möglichst durch empirische Theorien zu beschreiben, ist folgende Aussage Johnsons von besonderem Interesse:

\begin{small}
\enquote{..... I also want to anticipate the standard objection that, since we are bound to talk about preconceptual and nonpropositional aspects of experience always in propositional terms, it must follow that they are themselves propositional in nature. This simply doesn't follow. (...) ..... while we must use propositional language to describe these dimensions of experience and understanding, we must not mistake our mode of description for the things described.} [1987, p.4]   
\end{small}

Da wir nur die natürliche Sprache und formale Sprachen verfügbar haben, müssen wir uns mit deren Möglichkeiten begnügen, auch wenn wir damit nicht alle Aspekte unserer Erfahrung zum Ausdruck bringen können. In der Darstellung empirischer Theorien präzisieren die theoretischen Terme/Begriffe \enquote{preconceptual and nonpropositional aspects of experience and understanding}. Damit ändert sich nicht der Charakter dieser Aspekte, aber im Rahmen der Theorie werden sie kommunizierbar. Diese Vorgehensweise widerspricht auch keineswegs der Auffassung von Johnson, wie folgendes Zitat belegt:

\begin{small}
\enquote{Although we can more or less successfully abstract from particular empirical contents in framing a formal system (and in mathematics do it quite thoroughly), .....} [p. 38] \medskip
\end{small}
  
Geht man mit Sfard [2008] davon aus, daß sich Wissen i.w. in Diskursen (im weitesten Sinne) entwickelt und mitteilt --- was für die Unterweisung von Anfängern sicherlich zutrifft und dem Unterrichtsgespräch und auch der Lehrperson besondere Bedeutung beimißt, so erfolgt auch die Zuweisung von Bedeutungen in Diskursen. Sie unterscheidet zwei wesentliche Formen des Diskurses: \smallskip  

\begin{small}
\noindent \enquote{..... actual reality communication (auch \emph{\enquote{actual reality discourse}}, die Verf.) may be perceptually mediated by the objects that are being discussed, whereas in the \emph{virtual reality discourse} (kursiv durch die Verf.) perceptual mediation is scarce and is only possible with the help of what is understood as symbolic substitutes of objects under consideration. This description should not be read as a statement on an ontological status of the \enquote{realms} underlying the discourses. In introducing the metaphor of the two realities, I was psychologically, rather than philosophically, minded. That is, the distinction was drawn with an eye to the differing actions and experiences of the participants of the discourses rather than to ontological questions.} [2000, p. 39] \smallskip
\end{small}

Es entspricht dem Charakter empirischer Theorien, daß in ihnen \enquote{actual reality discourses} vorherrschen. \enquote{virtual reality discourses} ---  zu denen die mathematischen Diskurse zählen --- sind jedoch erforderlich, wenn theoretische Elemente/Begriffe thematisiert werden. \medskip 

Übereinstimmend mit dem Gesagten dienen bei der Formulierung empirischer Theorien \enquote{actual reality discourses} metaphorisch als Vorlage für \enquote{virtual reality discourses}. Beispiele solcher actual reality discourses sind: \smallskip

\begin{small}
--- Subtraction Is Taking smaller collections from larger collections to \\ \noindent \hspace*{8mm} form other collections.

--- Addition of a Given Quantity Is Taking Steps a given distance to the \\ \noindent \hspace*{8mm} right (or forward). \smallskip
\end{small}

Diese und weitere Beispiele findet man bei Lakoff und Núñez [1997], wobei wir uns auf solche Beispiele beschränken möchten, deren \enquote{virtual reality discourses} sich auf Elementarmathematik beziehen. \medskip 

Wie wird den Zeichen $\underset{\bar {S} \in \mathbb{S}'}{\bigcup} \bar{S}$ eine objektive Bedeutung zugewiesen? Hat man unendlich viele verschiedene Strecken, so liefert der endliche Fall keine Anhaltspunkte, ihn auf den unendlichen zu übertragen. Die Piagetschen Abstraktionsformen (\emph{empirische, pseudo--empirische, reflektive Abstraktion}) helfen hier nicht weiter. Die Terme $\underset{\bar{S} \in \mathbb{S}'}{\bigcup}\bar{S}$ können nur innerhalb eines
\enquote{virtual reality discourse} als Objekte bestimmt werden. \smallskip

Offenbar erfolgt eine Bedeutungszuweisung so, daß --- wie oben schon gesagt --- die Zeichen in den Modellen der Theorie T$_{\mathbb{R}^+}$ \enquote{geeignet} operationalisiert werden in \enquote{geeigneter} Weise mit den T$_{\mathbb{R}^+}$--vortheoretischen Begriffen verknüpft werden.  \medskip 

Felix Klein erwähnt in seiner \enquote{Elementarmathematik vom höheren Standpunkt} eine Konzeption der Bruchrechnung, in der der Autor eine \enquote{geeignete} Operationalisierung als eine \enquote{Verabredung} bezeichnet, für die es nur \enquote{Plausibilitätsgründe} gibt. Offensichtlich teilte Klein diese Interpretation. [1933, S. 32] \medskip   

Der Inhalt des erforderlichen Diskurses sind die folgenden Schritte. \medskip 

\begin{footnotesize}
\noindent \emph{Bemerkung}: Wir behalten die Schreibweise $\bar{S}, \bar{T}$, ... für die Strecken aus $\mathbb{S}$ bei, um deutlich zu machen, daß es sich um Elemente des Theorie--Elementes $\big ($T$_{\mathbb{R}^+} \big )$ handelt. \medskip 

Da die Strecken für Schüler geometrische Objekte sind, können sie die folgende Definition empirisch verifizieren. \medskip

$\bar{S}, \bar{T} \in \mathbb{S}$ \smallskip

$\bar{S} \leq \bar{T} \underset{\text{def}}{\leftrightarrow} \bar{S} \subseteq \bar{T}$  \medskip

\noindent Damit folgt: \smallskip

$\big ($T$_{\mathbb{R}^+}$:1$\big )  \; (\mathbb{S} \, ; \leq $) ist eine reflexive lineare Ordnung. \\ 
  
\noindent Nun die \emph{Modelle}: \medskip 

\begin{tabular}{l l l}
M$\big($T$_{\mathbb{R}^+} \big ) $ & = & $\big \langle \langle \mathbb{S}, \big \lbrace \underset{\bar {S} \in \mathbb{S}'}{\bigcup} \bar{S} \, \arrowvert \, \mathbb{S}' \subseteq \mathbb{S}\big \rbrace \rangle \big \rangle $ mit \\

& (i) & $ \langle \mathbb{S}, \big \lbrace \underset{\bar{S} \in \mathbb{S}'}{\bigcup} \bar{S} \, \arrowvert \, \mathbb{S}' \, \subseteq \mathbb{S} \big \rbrace \rangle \in $ M$_{p} \big($T$_{\mathbb{R}^+} \big )$ \\ 
&(ii)& $ \emptyset \neq \mathbb{S}' \, \subseteq \mathbb{S} \Rightarrow \underset{x}{\bigvee} \big ( x = \underset {\bar{S} \in \mathbb{S}'}{\bigcup} \bar{S} \wedge (x \in \mathbb{S} \, \vee \, x = \underset{\bar{S} \in \mathbb{S}}{\bigcup} \bar{S}) \big )$ \\ 
\end{tabular} \medskip

Da die Strecken als offene Intervalle in $\mathbb{Q}^+$ aufgefaßt werden können, erhält man \medskip

$\big ( $T$_{\mathbb{R}^+}$:2$\big ) \;  \mathbb{Q}^+ = \underset{\bar{S} \in \mathbb{S}}{\bigcup} \bar{S}$ 
\bigskip

\noindent M$\big($T$_{\mathbb{R}^+} \big)$(ii) läßt sich mit $\big ( $T$_{\mathbb{R}^+}$:2$\big )$ wie folgt formulieren: \smallskip

\begin{center}
$\emptyset \neq \mathbb{S}' \subseteq \mathbb{S} \Rightarrow \underset{\bar{S} \in \mathbb{S}'}{\bigcup} \bar{S}$ ist eine Strecke oder $\mathbb{Q}^+$ \\ 
\end{center}

\noindent oder 

\begin{center}
$ \emptyset \neq \mathbb{S}' \subseteq \mathbb{S} \, \wedge \underset{\bar{S} \in \mathbb{S}'}{\bigcup}\bar{S} \neq \mathbb{Q}^+ \Rightarrow \underset{\bar{T} \in \mathbb{S}}{\bigvee} \, \underset{\bar{S} \in \mathbb{S}'}{\bigwedge} (\bar{S} \leq \bar{T})$ \\
\end{center}

\sloppy Bedingung (ii) verknüpft folglich die theoretischen Elemente der Sprache mit den nichttheoretischen Anteilen. Damit ist die Adäquatheit von (ii) gegeben.
Da der Nachweis von $\big ( $T$_{\mathbb{R}^+}$:2$\big )$ keine Vereinigung von Strecken benötigt, wird durch die Umformulierung von M$\big($T$_{\mathbb{R}^+} \big)$(ii) das Problem, das Ergebnis unendliche vieler Vereinigungen zu bestimmen, auf elegante Weise umgangen. \\

\begin{tabular}{l l}
 Def.: & $ \emptyset \neq \mathbb{S}' \subseteq \mathbb{S} $ \smallskip \\
& Gilt $ \underset{\bar{S} \in \mathbb{S}'}{\bigcup} \bar {S} \neq \mathbb{Q}^+$, so heißt $\mathbb{S}' \,$ \emph{nach oben beschränkt} und jedes \\ 
& $\bar{T} \in \mathbb{S} $ mit $ \underset{\bar{S} \in \mathbb{S}'}{\bigcup} \bar{S} \leq \bar{T}$ heißt eine \emph{obere Schranke} zu $ \mathbb{S}'$ \bigskip \\ 
\end{tabular}

\begin{tabular}{l l}
$ \big ($T$_{\mathbb{R}^+}$:3$\big)$ & $ \emptyset \neq \mathbb{S}' \subseteq \mathbb{S} $ \smallskip \\
& Beh.: Ist $\mathbb{S}'$ nach oben beschränkt, so ist $\underset{\bar{S} \in \mathbb{S}'}{\bigcup} \bar{S}$ die \\  
& kleinste obere Schranke zu $\mathbb{S}'$,  das \emph{Supremum} von $\mathbb{S}'$. \smallskip \\
& bez.: sup $\mathbb{S}'$ \\
\end{tabular}

\begin{tabular}{l c r} 
M$\big ($T$_{\mathbb{R}^+} \big )$(ii) $\, \Rightarrow  \underset{\bar{S} \in \mathbb{S}'}{\bigcup}\bar{S}$ ist eine obere Schranke zu $\mathbb{S}'$ \\
\end{tabular}

\begin{tabular}{l c r}
\noindent $\bar{T} \in \mathbb{S} \wedge \underset{\bar{S} \in \mathbb{S}'}{\bigwedge} \big (\bar{S} \leq \bar{T} \big )$ & $   \underset{\text{def}}{\rightarrow} $ & $ \underset{\bar{S} \in \mathbb{S}'}{\bigwedge} \big ( \bar{S} \subseteq \bar{T} \big ) $ \\
& $ \Rightarrow $ & $ \underset{\bar{S} \in \mathbb{S}'}{\bigcup} \bar{S} \subseteq \bar{T} $ \\ 
& $ \underset{\text{def}}{\rightarrow} $ & $ \underset{\bar{S} \in \mathbb{S}'}{\bigcup} \bar{S} \leq \bar{T} $ \bigskip \\ 
\end{tabular}
\end{footnotesize}

Durch $ \big ($T$_{\mathbb{R}^+}$:3$\big)$ werden die $ \underset{\bar{S} \in \mathbb{S}'}{\bigcup} \bar{S}$ \enquote{geeignet} operationalisiert. \enquote{Geeignet} heißt, daß die Teilnehmer des Diskurses durch die Formulierung der Axiome sicherstellen, daß die Theorie die intendierten Anwendungen zutreffend beschreibt. Durch die Axiome erhalten die theoretischen Terme/Begriffe ihre Bedeutung innerhalb der Theorie, werden sie in die informelle Semantik der Theorie eingebunden. Man erwirbt sie \enquote{mit einem Schlag}, durch ein \enquote{Aha!}--Erlebnis, eine gelungene Anwendung. Die Modelle garantieren durch die Axiome, daß die Theorie erfolgreich angewendet werden kann. \smallskip

Die Konstruktion wird wie folgt abgeschlossen: Sei $\emptyset \neq \mathbb{S}' \subseteq \mathbb{S}$ und $\mathbb{S}'$ sei nach oben beschränkt. Es können zwei Fälle auftreten: \\
Im ersten Fall existiert ein kleinstes q $\in \mathbb{Q}^+$ so daß für jedes $\bar{S} \in \mathbb{S}' $ gilt $\bar {S} \subseteq$ (0,q) $\subseteq \mathbb{Q}^+$ Dann ist auch für den Schüler klar, daß jedes q$'$ < q und jedes q$''$ > q von sup $\mathbb{S}'$   verschieden ist. Trotzdem ist es für ihn (und auch für den Studenten der Anfängervorlesung) mit besonderen Schwierigkeiten verbunden, q mit sup $\mathbb{S}'$  gleichzusetzen. Denn sup $\mathbb{S}'$ wäre damit eine empirisch nicht verifizierbare Darstellung von q. Dies widerspricht seinem bisherigen empirisch geprägten Zahlverständnis. 

\noindent Im zweiten Fall, wenn es ein solches q $\in \mathbb{Q}^+$ nicht gibt, sup $\mathbb{S}'$ nicht der Ausgangsmenge $\mathbb{Q}^+$ angehört, sind die Verständnisschwierigkeiten noch größer. Denn auch in diesem Falle wird das Supremum als Zahl betrachtet --- man nennt sup $\mathbb{S}'$ eine \emph{irrationale Zahl}, als ein eigenständiges Objekt, völlig losgelöst von dem Prozeß, der es geschaffen hat (der Wechsel des \enquote{embodied schema}). Man kann dies in der Tat als einen ontologischen Wandel, als einen Wechsel im Seinsverständnis auffassen. \smallskip

Offen ist noch die Frage nach einer geeigneten Unterrichtskonzeption, in die sich die Vermittlung mathematischen Wissens im Rahmen empirischer Theorien einordnen läßt.

Als in den 1960er Jahren die Neuorientierung der Anwendungen (der \enquote{eingekleideten Aufgaben}) diskutiert wurde, bestimmte der Begriff des \emph{Mathematisierens} die Diskussion. In der französischen Mathematikdidaktik sprach man von \emph{mathématisation des situations}. Im Gegensatz zur Modellbildung, die auf der Trennung von Realität und Mathematik fußt, wurde bei der mathématisation des situations das Ineinandergreifen von Realität und Mathematik betont:   

\begin{small}
\enquote{La mathématique est faite par les abstractions, mais elle change la réalité physique, elle s'y intègre. Elle n'est pas à c$\hat{o}$té du réel, elle est dans le réell, elle lui donne une nouvelle dimension. ..... Il n'y a pas de fossé entre le concret et l'abstrait. Ces deux mots détermint deux positions méthodologiques différentes de l'homme par rapport à la réalité quand il l'explore et la change. \\ 
Il y a une dialectique intrinsèque dans la recherche de la réalité. Nous faisons les abstractions pour mieux conna$\hat{i}$tre le concret et pour le changer. De cette fa\c{c}on, nous enrichissons la réalité. Le concret de cette réalité exige de nouvelles abstractions pour $\hat{e}$tre mieux compris et changé. Et ainsi de suite.} (Smolec [1972, p. 125])  
\end{small}

Auf der ICMI Tagung \emph{How to Teach Mathematics so as to Be Useful} in Utrecht 1967 formuliert Anna Krygowska:  

\begin{small}
\enquote{Il arrive que certaines constructions mathématiques aient des réalisations matérielles particulièrement directes dans ce sens qu'elles puissent servir directment à la description de la réalité et à la prise de décisions concrètes.} [1968/69, p. 10]  
\end{small}

\noindent Freudenthal bezieht diesen Gedanken auf die Arithmetik:  

\begin{small}
\enquote{ ..... arithmetic starts in a concrete context and patiently it returns to concrete contexts as often as needed.} [1968/69, p. 6] 
\end{small}

\noindent Von Griesel wird er wie folgt formuliert:  

\begin{small}
\enquote{Die mathematischen Gegenstände sind nach meiner Meinung vom einzelnen Menschen gebildete, also im einzelnen Menschen durch Gehirnvorgänge verkörperte, fiktive, gedankliche Konstrukte, Konzepte. ..... \smallskip \\
Diese Konstrukte werden wenigstens im Bereich der Elementarmathematik in aktiver Auseinandersetzung mit der Erfahrungswirklichkeit in einem Akt \emph{ideeisierender Abstraktion} gebildet und können dann wieder in diese \emph{hineininterpretiert} (der Welt \enquote{übergestülpt}) werden. \smallskip \\
Diese Fiktionen haben dann \emph{Teil an der Erfahrungswirklichkeit}. Sie sind nach der Hineininterpretation in diese \emph{eingebettet. Sie sind real inhärent in der Erfahrungswirklichkeit}.} [2013b]  \smallskip 
\end{small}

Auch in der Didaktik der früheren Volksschule (weitgehend übereinstimmend mit der Grundschule und der heutigen Sekundarstufe I nichtgymnasialer Schulformen; die Verf.) waren diese Gedanken lebendig. Bei Wilhelm Oehl liest man: \smallskip 

\begin{small}
\enquote{Damit ein Kind einen mathematischen Begriff auf Sachsituationen anwenden kann, muß es vorher durch wiederholtes eigenes Tun erlebt haben, wie dieser Begriff aus Sachsituationen von ganz bestimmter Struktur herauswächst. Diese konkret anschaulichen Handlungen allein können der Altersstufe des Kindes gemäß die Grundlage neu zu gewinnender mathematischer Begriffe sein. So kommen wir zu der Forderung: \so{Im Rechenunterricht der Volksschule mu\ss{ } jeder neue Begriff aus der t\"atigen Anschauung heraus entwickelt werden.} .....  \\
Erst wenn die neuen Begriffe auf konkret anschaulichem Weg gewonnen sind, werden ihre Beziehungen zu den bereits bekannten Begriffen hergestellt. ..... \\

Kennzeichnend für diese Art des didaktischen Vorgehens im Rechenunterricht der Volksschule sind zwei Kriterien: \so{Die Situationsgebundenheit des rechnerischen Tuns} und die \so{Sachgebundenheit des mathematischen Denkens.} Unter der Situationsgebundenheit ist der natürliche Spannungsbogen zwischen Kind und Welt zu verstehen, der den Charakter eines \enquote{wirkenden Tuns}\footnote{F. Hartke, Die psychologische Fragestellung. Pädag. Rundschau, H. 7 / 1961, S. 403} annimmt, eines auf den Sachgegenstand bezogenen Verhaltens, das dann gegeben ist, wenn der Schüler die Sachsituation durchschaut und meistert. Dieser tätige Umgang mit den Dingen, der die Voraussetzung für das rechnerische Tun darstellt, verlagert sich im didaktischen Vollzug von der Ebene des ursprünglichen, natürlichen Verhaltens auf die des klärenden Überlegens. Die Dinge mit \enquote{Umgangsqualitäten} werden dann zu Sachgegenständen für mathematische Einsicht: Das ist die sich aus der Situation entwickelnde \so{Sachgebundenheit} des \so{mathematischen Denkens}.} [1962, S. 16] \medskip    
\end{small}

Diese  Gedanken entsprechen genau der Situation, die Kitcher skizziert. Der Schüler beschreibt nicht nur mit mathematischen Mitteln die Realität, sondern diese befruchtet und erweitert gleichzeitig sein mathematisches Wissen, z.B. dadurch, daß die Gegebenheiten, mit denen er sich auseinandersetzt, leicht variieren. Wie man das Addieren von 2 und 3 durch Zusammenschieben und Abzählen von Steinen erlernen kann so auch durch Zusammenschieben und Abzählen von Perlen das Addieren von 2 und 4. \enquote{Knowing is doing} lautet ein Aphorismus von Humberto Maturana und Francisco Varela (Brent Davis [1995, p. 4]). Das Wechselspiel zwischen Realität und Mathematik ist bei diesen Lernprozessen wesentlich intimer als bei einer expliziten Modellbildung. Es betont die Verzahnung von Realität und Mathematik und nicht ihre Trennung. Ein in diesem Sinne verstandenes Mathematisieren entwickelt somit genau das Wissen, dessen Inhalte in den empirischen Theorien erworben wird, auf die sich die Elementarmathematik stützt. Und es entspricht der oben skizzierten Auffassung, Mathematik wie eine Naturwissenschaft zu betrachten. Das Wissen, das die Kinder gewinnen, entwickelt sich im Rahmen empirischen Vorgehens, basierend auf Daten und Entdeckungen. 

Diese Auffassung deckt sich mit der Sichtweise der Historiker. Mit den Worten von Gerhard Frey:  

\begin{small}
\enquote{Mathematik ist nicht als reine Mathematik entstanden, sondern in und gleichzeitig mit ihren Anwendungen auf die Wirklichkeit.} [1967, S. 7] 
\end{small}

\vspace{2cm}

\selectlanguage{german}

\nocite{*}

\printbibliography

\end{document}